Title:  Almost Complex Structures on $S^2\times S^2$
Author: Dusa McDuff
Subj-class: SG
MSC-class: 53C15
Comments: Latex 44 pages, no pictures
Abstract:
In this note we investigate the structure of the space $\Jj$ of
smooth almost complex structures on $S^2\times S^2$ that are compatible
with  some symplectic form.  This space has a natural stratification that
changes as the cohomology class of the form changes and whose properties
are very closely connected to the topology of the group of
symplectomorphisms of $S^2\times S^2$.
  By globalizing standard gluing
constructions in the theory of  stable maps, we show that the strata of $\Jj$ 
 are Fr\'echet manifolds of finite
codimension, and that the normal link of each  stratum is a  finite
dimensional stratified space.   The topology of these links turns out to be
surprisingly intricate, and we work out  certain cases. 
Our arguments apply also to other ruled surfaces, though they give
complete information only for  bundles over $S^2$ and $T^2$.

\documentstyle{article}

\newcommand{\R}{{\bf R}}
\newcommand{\ov}{\overline}
\newcommand{\bcp}{{\overline{{\C}P}\,\!^{2}}}
\newcommand{\pl}{{\mbox{$\,\equiv \hspace{-.92em} \| \,\,\,$}}}
\newcommand{\p}{{\partial}}
\newcommand{\oMm}{{\ov{\Mm}}}

\newcommand{\op}{{\overline{\partial}}}
\newcommand{\Coker}{{\rm Coker}}
\newcommand{\Ker}{{\rm Ker}}
\newcommand{\id}{{\rm Id}}
\newcommand{\PSL}{{\rm PSL}}
\newcommand{\oJj}{{\ov{\cal J}}}
\newcommand{\ind}{{\rm ind\,}}

\newcommand{\TL}{{\Tilde L}}
\newcommand{\TY}{{\Tilde Y}}
\newcommand{\TU}{{\Tilde U}}
\newcommand{\Tilde}{\widetilde}
\newcommand{\tsi}{{\tilde{\si}}}

\newcommand{\Aut}{{\rm Aut}}

\newcommand{\G}{{\rm G}}

\newcommand{\im}{{\rm Im\,}}
\newcommand{\Symp}{{\rm Symp}}

\newcommand{\Ee}{{\cal E}}
\newcommand{\Kk}{{\cal K}}
\newcommand{\Mm}{{\cal M}}
\newcommand{\Nn}{{\cal N}}
\newcommand{\Ll}{{\cal L}}
\newcommand{\Dd}{{\cal D}}
\newcommand{\Xx}{{\cal X}}
\newcommand{\Jj}{{\cal J}}
\newcommand{\Cc}{{\cal C}}
\newcommand{\Ff}{{\cal F}}
\newcommand{\Gg}{{\cal G}}
\newcommand{\TGg}{{\Tilde {\Gg}}}
\newcommand{\Uu}{{\cal U}}

\newcommand{\Zz}{{\cal Z}}
\newcommand{\Vv}{{\cal V}}
\newcommand{\Pp}{{\cal P}}

\newcommand{\Ss}{{\cal S}}
\newcommand{\Oo}{{\cal O}}
\newcommand{\th}{{\theta}}

\newcommand{\Z}{{\bf Z}}
\newcommand{\C}{{\bf C}}

\newcommand{\al}{{\alpha}}

\newcommand{\be}{{\beta}}

\newcommand{\om}{{\omega}}
\newcommand{\eps}{{\varepsilon}}
\newcommand{\de}{{\delta}}
\newcommand{\De}{{\Delta}}
\newcommand{\ga}{{\gamma}}
\newcommand{\Ga}{{\Gamma}}
\newcommand{\ka}{{\kappa}}
\newcommand{\la}{{\lambda}}
\newcommand{\si}{{\sigma}}
\newcommand{\La}{{\Lambda}}

\newcommand{\io}{{\iota}}
\newcommand{\Si}{{\Sigma}}

\newcommand{\INT}{{\rm Int}}

\newcommand{\SS}{{\smallskip}}
\newcommand{\MS}{{\medskip}}

\newcommand{\NI}{{\noindent}}
\newcommand{\proof}[1]{\noindent{\bf Proof#1:\  }}

\newcommand{\QED}{\hfill$\Box$\medskip}

\newtheorem{theorem}{Theorem}[section]
\newtheorem{thm}[theorem]{Theorem}

\newtheorem{defn}[theorem]{Definition}
\newtheorem{example}[theorem]{Example}
\newtheorem{remark}[theorem]{Remark}
\newtheorem{lemma}[theorem]{Lemma}

\newtheorem{prop}[theorem]{Proposition}

\begin{document}
\title{Almost complex structures on $S^2\times S^2$}
\author{Dusa McDuff\thanks{Partially
supported by NSF grant DMS 9704825.} \\ State University of New York
at Stony Brook \\ (dusa@math.sunysb.edu)}

\date{Aug 1, 1998}
\maketitle
\MS

\MS

\NI
{\bf Abstract}

In this note we investigate the structure of the space $\Jj$ of
smooth almost complex structures on $S^2\times S^2$ that are compatible
with  some symplectic form.  This space has a natural stratification that
changes as the cohomology class of the form changes and whose properties
are very closely connected to the topology of the group of
symplectomorphisms of $S^2\times S^2$.
  By globalizing standard gluing
constructions in the theory of  stable maps, we show that the strata of $\Jj$ 
 are Fr\'echet manifolds of finite
codimension, and that the normal link of each  stratum is a  finite
dimensional stratified space.   The topology of these links turns out to be
surprisingly intricate, and we work out  certain cases. 
Our arguments apply also to other ruled surfaces, though they give
complete information only for  bundles over $S^2$ and $T^2$.

\section{Introduction}

It is well known that every symplectic form on $X = S^2\times S^2$ is, after
multiplication by a suitable constant,  symplectomorphic to a product form
$\om^\la = (1+\la)\si_1 + \si_2$ for some $\la \ge 0$, where the  $2$-form
$\si_i$ has total area $1$ on the $i$th factor.  We are interested in the
structure of the space $\Jj^\la$ of all $C^\infty$ $\om^\la$-compatible almost
complex structures on $X$.  Observe that $\Jj^\la$ itself is always
contractible.  However it has a natural stratification that
changes as $\la$ passes each integer.  The reason for this is that as $\la$
grows the set of homology classes that can be represented by an
$\om^\la$-symplectically embedded $2$-sphere changes.  Since each such
$2$-sphere can be parametrized to be $J$-holomorphic for some $J\in
\Jj^\la$, there is a corresponding change in the structure of $\Jj^\la$.

To explain this in more detail, let $A\in H_2(X,\Z)$ be the homology class
$[S^2\times pt]$ and let $F = [pt\times S^2]$.  (The reason for this notation is
that we are thinking of $X$ as a fibered space over the first $S^2$-factor, so
that the smaller sphere $F$ is  the fiber.) When $\ell -1 < \la \le \ell$,
$$
\om^\la(A-kF) > 0, \quad\mbox{for}\;\; 0\le k \le \ell.
$$
Moreover, it is not hard to see that for each such $k$ there is a map
$\rho_k:S^2\to S^2$ of degree $-k$ whose  graph
$$
z\mapsto (z,\rho_k(z))
$$
is an $\om^\la$-symplectically embedded sphere in $X$.  It follows easily that
the space 
$$
\Jj_k^\la = \{J\in \Jj^\la: \mbox{there is a $J$-hol curve in class \,} A-kF\}
$$
is nonempty  whenever $k< \la + 1$.  Let
$$
\oJj_k^\la = \cup_{m\ge k} \;\Jj_m^\la.
$$
Because $(A-kF)\cdot (A-mF) < 0$ when $k\ne m > 0$, positivity of
intersections implies that there is exactly one $J$-holomorphic curve in class
$A-kF$ for each $J\in \Jj_k$.  We denote this curve by $\Delta_J$.

\begin{lemma}\label{basic} The spaces $\Jj_k^\la, 0\le k \le \ell$, are
disjoint  and $\oJj_k^\la$ is the closure of $ \Jj_k^\la$ in $\Jj^\la$.  Further,
$\Jj^\la = \oJj_0^\la$. 
\end{lemma} 
\proof{} It is well known that for every $J\in \Jj$ the
set of $J$-holomorphic curves in class $F$ form the fibers of a fibration $\pi_J:
X\to S^2$.  Moreover, the class $A$ is represented by either a curve or a
cusp-curve (i.e. a stable map).\footnote
{
We will follow the convention of [LM]  by defining a ``curve" to be  the
image of a single sphere, while a ``cusp-curve" is either
multiply-covered or has domain equal to a union of two or more spheres.}
Since the class $F$ is always represented and $(mA + pF)\cdot F = m$, it
follows from  positivity of intersections that $m\ge 0$ whenever $mA+pF$
is represented by a curve. Hence any
  cusp-curve in class $A$  has one component in some class
$A-kF$ for $k\ge 0$,  and  all others represent a multiple of $F$. 
In particular, each $J\in \Jj^\la$ belongs to some set $\Jj_k^\la$.  Moreover,
because $(A-kF)\cdot (A-mF) < 0$ when $k\ne m$ and $k,m\ge 0$, the
different $\Jj_k^\la$ are disjoint. 

The second statement holds because if $J_n$ is a sequence of elements 
in $\Jj_k^\la$, then the corresponding sequence of $J_n$-holomorphic
curves in class $A-kF$ has a convergent subsequence 
whose limit is a cusp-curve in class $A-kF$.  This limit has to have a
component in some class $A-mF$, for $m\ge k$, and so $J\in \Jj_m^\la$ for
some $m\ge k$.   For further details see Lalonde--McDuff [LM], for
example.\QED

Here is our main result.  Throughout we are working with
 $C^\infty$-maps and almost complex structures, and so by manifold we mean a
Fr\'echet manifold.  By a stratified space $\Xx$ we mean a topological space
that is a union of a finite number of disjoint manifolds that are called strata. 
Each stratum $\Ss$ has a neighborhood $\Nn_\Ss$ that projects to $\Ss$ by a
map $\Nn_\Ss\to \Ss$.  When $\Nn_\Ss$ is given the
induced stratification, this map is a locally trivial
fiber bundle whose fiber has the
form of a cone $C(\Ll)$  over a finite dimensional stratified space $\Ll$ that is
called the {\em link} of $\Ss$ in  $\Xx$.   Moreover, $\Ss$ sits inside $\Nn_\Ss$
as the set of vertices of all these cones.

\begin{thm}\label{main} (i)  For each $1\le k\le \ell$, $\Jj_k^\la$ is a
submanifold of $\Jj^\la$ of codimension $4k-2$.
\SS

\NI
(ii)  For each $m > k \ge 1$ the normal link $\Ll_{m,k}^\la$ of $\Jj_{m}^\la$ in
$\oJj_k^\la$ is a stratified space of dimension $4(m-k) -1$.  Thus,
there is a neighborhood of $\Jj_{m}^\la$ in $\oJj_k^\la$ that is fibered over
$\Jj_{m}^\la$ with fiber equal to the cone on $\Ll_{m, k}^\la$.
\SS

\NI
(iii)  The structure of the link $\Ll_{m,k}^\la$  is independent of $\la$ (provided
that $\la >m -1$.)  
 \end{thm}

The first part of this theorem was proved by Abreu in~[A], at least in the
$C^s$-case where $s<\infty$.  (Details are given in \S4.1 below.)  The second
and third parts follow by globalising recent work by Fukaya--Ono [FO],
Li--Tian~[LiT], Liu--Tian~[LiuT1], Ruan~[R] and others  on the structure of the
compactification of moduli spaces of $J$-holomorphic spheres via stable
maps.  We extend current gluing methods by showing that it is possible to
deal with obstruction bundles whose elements do not vanish at the gluing
point: see \S 4.2.3.  Another essential point is that we use Fukaya and
Ono's method  of dealing with the ambiguity in the parametrization of a
stable map since this involves the least number of choices and 
allows us to   globalize by constructing a gluing map that is equivariant with
respect to suitable local torus actions: see~\S 4.2.4 and \S 4.2.5.

The above theorem is
the main tool used in [AM] to calculate the rational cohomology ring of the
group $G^\la$ of symplectomorphisms of $(X, \om^\la)$.  

Observe that  part (iii)  states that the normal
structure of the stratum $\Jj_k^\la$ does not change with $\la$.  On the other
hand, it follows from the results of [AM] that the cohomology  of $\Jj_k^\la$
definitely does change as $\la$ passes each integer.  Obviously, it would be
interesting to know if the topology of $\Jj_k^\la$ is otherwise fixed. For
example, one could try to construct maps $\Jj^\la\to \Jj^{\mu}$ for $\la < \mu$
that preserve the stratification, and  then try to prove that they induce
homotopy equivalences $\Jj_k^\la \to \Jj_k^{\mu}$ whenever $\ell -1 < \la \le
\mu \le \ell$.   

The most we have so far
managed to do in this direction is to prove the following lemma that, in
essence, constructs  maps  $\Jj^\la\to \Jj^{\mu}$ for $\la < \mu$.
It is not clear whether these are homotopy equivalences for $\la, \mu \in
(\ell-1, \ell]$. It is convenient to fix a
fiber $F_0 = pt\times S^2$ and define 
$$
\Jj_k^\la(\Nn(F_0)) = \{ J\in \Jj_k^\la: J = J_{split}\mbox{ near } F_0\},
$$
where $J_{split} $ is the standard product almost complex structure.

\begin{lemma}\label{inc} (i) The inclusion $\Jj^\la(\Nn(F_0))\to \Jj^\la$
induces a homotopy equivalence $\Jj_k^\la(\Nn(F_0))\stackrel{\simeq}{\to}
\Jj_k^\la$ for all $k <  \la+1$.

\NI
(ii)  Given any compact subset $C\subset \Jj^\la(\Nn(F_0))$ and any
$ \mu> \la$ there is a map
$$
\io_{\la,\mu}:  C\to \Jj^{\mu}(\Nn(F_0))
$$
that takes $C\cap \Jj_k^\la(\Nn(F_0))$ into $\Jj_k^{\mu}(\Nn(F_0))$ 
for all $k$.
\end{lemma}

This lemma is proved in \S 2.

 The next task is to calculate the
links $\Ll_{m,k}^\la$.  So far this has been done for the easiest case:

\begin{prop}\label{pr:lens}  For each $k\ge 1$ and $\la$ the link 
$\Ll_{k+1,k}^\la$ is
the $3$-dimensional lens space $L(2k,1)$. 
\end{prop}

Finally, we illustrate our methods by using the stable map approach to
confirm that the link of $\Jj_2^\la$ in $\Jj^\la$ is $S^5$,
 as predicted by part (i) of
Theorem~\ref{main}.   Our method first calculates an auxiliary link $\Ll_\Zz$
from which the desired link is obtained by collapsing certain strata.  
The $S^5$ appears in a  surprisingly interesting way, that can be  briefly
described as follows.

Let $\Oo(k)$
denote the complex line bundle over $S^2$ with Euler number $k$, where we
 write $\C$ instead of $\Oo(0)$.  Given a vector bundle $E\to B$ we write
$S(E)\to B$ for its unit sphere bundle.  Note that
 the unit $3$-sphere bundle 
$$
S(\Oo(k) \oplus \Oo(m))\to S^2
$$
 decomposes as the composite 
$$
S(L_{P(k,m)})\to \Pp (\Oo(k)\oplus \Oo(m)) \to  S^2
$$
where  $L_{P(k,m)}\to \Pp (\Oo(k)\oplus \Oo(m))$ is the canonical line bundle 
over the the projectivization  of $\Oo(k)\oplus \Oo(m)$.  
In particular,  the space
$S(\Oo(-1) \oplus \C)$ can be identified with $S(L_{P(-1,0)})$.
But $\Pp (\Oo(-1)\oplus \C)$ is simply the blow up $\C P^2\#\bcp$,
and its canonical bundle is the pullback of the canonical bundle over
$\C P^2$.   We also consider the  singular line bundle 
(or orbibundle) $L_Y\to L$ whose associated unit sphere bundle has total space
$L(Y) = S^5$ and fibers equal to the orbits of the following $S^1$-action on
$S^5$:    $$  \th\cdot (x,y,z) = (e^{i\th}x, e^{i\th}y, e^{2i\th}z),\quad x,y,z\in \C. 
$$

\begin{thm}\label{LINK} (i)
The space $\Ll_\Zz$ obtained by plumbing the unit sphere bundle of
$\Oo(-3)\oplus \Oo(-1)$ with the singular circle bundle $S(L_Y)\to Y$ may be
identified with the unit circle bundle of the canonical bundle
over $\Pp(\Oo(-1)\oplus \C) = \C P^2\#\bcp$.  

\NI
(ii)  The link $\Ll_{2,0}^\la$ is
obtained from $\Ll_\Zz$ by collapsing the fibers  over the exceptional divisor
 to a single fiber, and hence may be identified with
$S^5$.  Under this identification, the link $\Ll_{2,1}^\la = \R P^3$ corresponds
to the inverse image of a conic in $\C P^2$. \end{thm}

In his recent paper~[K], Kronheimer shows that the universal deformation of
the quotient singularity $\C^2/(\Z/m\Z)$ is transverse to all the submanifolds
$\Jj_k$ and so is an explicit model for the normal slice of $\Jj_m$ in $\Jj$.
  Hence one can investigate
the structure of the intermediate links $\Ll_{m,k}^\la$ using tools from
algebraic geometry.  It is very possible that it would be  easier to
calculate these links this way.  However, it is still  interesting to try to
understand these links from the point of view of stable maps, since this is more
closely connected to  the symplectic geometry of the manifold $X$.

  Another point is that throughout we consider $\om^\la$-compatible
almost complex structures rather than $\om^\la$-tame ones.  However, it is
easy to see that all our results hold in the latter case.  
\MS

\subsection*{Other ruled surfaces}

All the above results have analogs for other ruled surfaces $Y\to \Si$.
If $Y$ is diffeomorphic to the product $\Si\times S^2$, we can define
$\om^\la, \Jj_k^\la$ as above, though now we should allow $\la$ to be any
number $> -1$ since there is no symmetry between the class $A = [\Si\times
pt]   $ and $F = [pt\times S^2]$.  In this case Theorem~\ref{main} still holds.
The reason for this is that if $u:\Si\to Y$ is an injective  $J$-holomorphic map in
class $A-kF$ where $k\ge 1$, then the normal bundle $E$ to the image
$u(\Si)$ has negative first Chern class so that the linearization $Du$ of $u$
has kernel and cokernel of constant dimension.  (In fact, the normal part of
$Du$ with image in $E$ is injective in this case.  See Theorem 1$^\prime$ in
Hofer--Lizan--Sikorav~[HLS].)   However, Lemma~\ref{basic} fails unless
$\Si$ is a torus since  there are tame almost complex structures on $Y$ with no
curve in class $[A]$.  One might think to remedy this by adding other strata
$\Jj_{-k}^\la$ consisting of all $J$ such that the class $A+kF$ is
represented by a $J$-holomorphic curve $u:(\Si, j)\to Y$ for some complex
structure $j$ on $\Si$.   However, although the universal moduli space
$\Mm(A+kF, \Jj^\la)$ of all such pairs $(u,J)$ is a  manifold, the map $(u,J)\to J$
is no longer injective: even if one cuts down the dimension by fixing a suitable
number of points each $J$ will in general admit several curves through these
points.   Moreover,  as $u$ varies over $\Mm(A+kF, \Jj^\la)$ the dimension of
the kernel and cokernel of $Du$ can jump.  Hence the argument given in 
\S 4.1 below that the strata $\Jj_k^\la$ are submanifolds of
$\Jj^\la$ fails on several counts.

In the case of the torus, $\Jj_0^\la$ is open and so Lemma~\ref{basic} does
hold.  However, it is not clear whether this is enough for the main application,
which is to further our
understanding of the groups $\G^\la$ of symplectomorphisms of $(Y, \om^\la)$. 
One crucial ingredient of the argument in~[AM] is that the action of this
group  on each stratum $\Jj_k^\la$ is essentially transitive.  More precisely, we
show that the action of $\G^\la$ on $\Jj_k^\la$ induces a homotopy equivalence
 $\G^\la/\Aut(J_k) \to
\Jj_k^\la$, where $J_k$ is an integrable element of
$\Jj_k^\la$ and $\Aut(J_k)$ is its stabilizer.  It is not clear whether this would
hold for the stratum $\Jj_0^\la$ when $\Si = T^2$.  One might have to take into
account the finer stratification considered by Lorek in~[Lo].  
He points out that the space
$\Jj_0^\la$   of all $J$ that admit a curve in class $A$ is not homogeneous.  A
generic element admits a finite number of such curves that are regular (that is
$Du$ is surjective), but since this number can vary the set of regular elements
in $\Jj_0^\la$ has an infinite number of components.  Lorek also characterises
the other strata that occur.  For example, the codimension $1$ stratum consists
of $J$ such that all $J$-holomorphic $A$ curves are isolated but there is at least
one where the kernel of $Du$ has dimension $3$ instead of $2$.  (Note that 
these $2$ dimensions correspond to the reparametrization group, since $Du$ is
the full linearization, not just the normal component.)

Similar remarks can be made about the case when $Y\to \Si$ is a
nontrivial bundle.  In this case we can label the strata $\Jj_k^\la$ so that the
$J\in \Jj_k^\la$ admit sections with self-intersection $-2k+1$.   Again 
Theorem~\ref{main} holds, but Lemma~\ref{basic} may not.  When $\Si
= S^2$ the homology class of the exceptional divisor is always represented, so
that $\Jj^\la = \oJj_1^\la$.  When $\Si = T^2$, the homology class of the section
of self-intersection $+1$ is always represented.   Thus $\Jj^\la = \oJj_{-1}^\la$.
Hence the analog of Lemma~\ref{basic} holds in these two
cases.  Moreover all embedded tori of self-intersection $+1$ are 
regular (by the same result in~[HLS]), which may help in the application
to $\Symp(Y)$.

We now state in detail the result for the nontrivial bundle $Y\to S^2$ since this
is used in~[AM].  Here  $Y= \C P^2\#\bcp$, and so every symplectic form on
$Y$ can be obtained from an annulus $A_{r,s} =\{ z\in \C^2: r \le |z| \le s\}$ by
collapsing the boundary spheres to $S^2$ along the characteristic orbits. This
gives rise to a form $\om_{r,s}$ that takes the value $\pi s^2$ on the class $L$
of a line and $\pi r^2$ on the exceptional divisor $E$.  Let us write $\om^\la$
for the form      $\om_{r,s}$ where $\pi s^2 = 1+\la, \pi
r^2  = \la > 0$.     Then the class $F = L-E$ of the fiber has size
$1$  as before, and $\Jj_k^\la$, $k\ge 1$, is the set of
$\om^\la$-compatible $J$ for which the class $E - (k-1) F$ is
represented.  

\begin{thm} When $Y =  \C P^2\#\bcp$ the spaces
$\Jj_k^\la$ are Fr\'echet submanifolds of $\Jj^\la$ of codimension
$4k$, and form the strata of a stratification of $\Jj^\la$ whose normal
structure is independent of $\la$.  Moreover, the normal link of
$\Jj_{k+1}^\la$ in $\Jj_k^\la$ is the lens space $L(4k+1, 1)$, $k\ge 1$.
\end{thm}
\MS

This paper is organised as follows.  \S 2 describes the main ideas in the
proof of Theorem~\ref{main}.  This relies heavily on the theory of stable
maps, and for the convenience of the reader we outline its main points.
References for the basic theory are for example [FO], [LiT] and [LiuT1].
\S 3 contains a detailed calculation of the link of $\Jj_2^\la$ in $\Jj^\la$.
In particular we discuss  the topological structure of the space of
degree $2$ holomorphic self-maps of $S^2$ with up to $2$ marked
points, and of the canonical line bundle that it carries. Plumbing with the
 orbibundle $L_Y\to Y$ turns out to be a kind of orbifold
blowing up process: see \S 3.1.  Finally, in \S 4 we work out the
 technical details of gluing that are needed to 
establish that the submanifolds $\Jj_k^\la$ do have 
a good normal structure.  The basic method here is
taken from McDuff--Salamon [MS] and Fukaya--Ono [FO].  
\MS

 \NI
{\bf Acknowledgements}\SS

I wish to thank Dan Freed, Eleni Ionel and particularly John Milnor for
useful discussions on various aspects of the calculation in \S 3, and Fukaya and
Ono for explaining to me various details of their arguments.

 \section{Main ideas}

We begin by proving Lemma~\ref{inc} since this is elementary,
and then will
 describe the main points in the proof of
Proposition~\ref{main}.  
\MS

\subsection{The effect of increasing $\la$}

\NI
{\bf Proof of Lemma~\ref{inc}}

Recall that  $F_0$ is a fixed fiber $pt\times S^2$ and that
$$
\Jj_k^\la(\Nn(F_0)) = \{ J\in \Jj_k^\la: J = J_{split}\mbox{ near } F_0\},
$$
We will also use the space
$$
\Jj_k^\la(F_0) = \{ J\in \Jj_k^\la: J = J_{split}\mbox{ on } TF_0\}.
$$
Let $\Ff^\la$ be the space of $\om^\la$-symplectically embedded curves in
the class $F$
 through a fixed point $x_0$.  Because there is a unique
$J$-holomorphic $F$-curve through $x_0$ for each $J\in \Jj$ (see
Lemma~\ref{basic}), there is a fibration
$$
\Jj^\la(F_0)\to \Jj^\la \to \Ff^\la.
$$
Since the elements of $\Jj^\la(F_0)$ are sections of a bundle with
contractible fibers, $\Jj^\la(F_0)$ is contractible.  Hence $\Ff^\la$ is also
contractible.  By using the methods of Abreu [A], it is not hard to show
that the symplectomorphism group $\Gg^\la = \Symp_0(X,\om^\la)$  of
$(X,\om^\la)$ acts transitively on $\Ff^\la$. Since the action of $\Gg^\la$
on $\Jj^\la$ preserves the strata $\Jj_k^\la$, it follows that the projection
$\Jj_k^\la\to \Ff^\la$ is surjective.  Hence there are induced fibrations $$
\Jj_k^\la(F_0)\to \Jj_k^\la \to \Ff^\la.
$$
This implies that the inclusion $\Jj_k^\la(F_0)\to \Jj_k^\la$ is a weak homotopy
equivalence.  

We now claim that the inclusion $\Jj_k^\la(\Nn(F_0))\to \Jj_k^\la(F_0)$ is also a
weak homotopy equivalence.  To prove this, we need to show that the
elements of any compact set $\Kk\subset \Jj_k^\la$ can be homotoped
near $F_0$ to make them coincide with $J_{split}$.  Since the set of tame
almost complex structures at a point is contractible, this is always possible
in $\Jj^\la$: the difficulty here is to ensure that $\Kk$ remains in
$\Jj_k^\la$ throughout the homotopy.  Here is a sketch of one method.
For each  $J\in \Jj_k^\la$ let $\De_J$ denote the unique
$J$-holomorphic curve in class $A-kF$.  Then $\De_J$ meets $F_0$
transversally at one point, call it $q_J$.  For each $J\in \Kk$, isotop the
curve $\De_J$ fixing $q_J$ to make it coincide 
 in a small neighborhood of $q_J$ with the flat section $S^2\times
pt$ that contains $q_J$.  (Details of an
very similar construction can be found in [MP], Prop 4.1.C.)  Now lift 
this isotopy  to  $\Jj_k^\la$.  Finally adjust the family of almost complex
structures near $F_0$, keeping $\De_J$ holomorphic throughout.

This proves (i).  Statement (ii) is now easy.  For any compact subset $C$
of $\Jj^\la(\Nn(F_0))$ there is $\eps>0$ such that $J = J_{split}$ on the
$\eps$-neighborhood $\Nn_\eps(F_0)$  of $F_0$.  
Let $\rho$ be a  nonegative
$2$-form supported inside the $2$-disc of radius $\eps$ that vanishes
near $0$, and let $\pi^*(\rho)$ denote its pullback to $\Nn_\eps(F_0)$ by the
obvious projection.  Then every $J$ that equals $J_{split}$ on
$\Nn_\eps(F_0)$ is compatible with  the form  $\om^\la + \ka\pi^*(\rho)$ for
all $\ka > 0$.   Since  $\om^\la + \ka\pi^*(\rho)$ is isotopic to $\om^{\mu}$ for
some $\mu$,  there is a diffeomorphism $\phi$ of $X$ that is isotopic to the
identity and is such that $\phi^*(\om^\la + \ka\pi^*(\rho)) = \om^{\mu}$. 
Moreover, because, by construction, $\pi^*(\rho) = 0$ near $F_0$,  we can
choose $\phi = \id $ near $F_0$.    Hence the map $ J\mapsto \phi^*(J) $
takes $\Jj^\la(\Nn(F_0))$ to $\Jj^{\mu}(\Nn(F_0))$.  Clearly it preserves
the strata $\Jj_k$.\QED

\subsection{Stable maps}

 From now on,
we will drop $\la$ from the notation, assuming that $k<\la + 1$ as before.  
We study  the spaces $\Jj_k$ and $\oJj_k$ by
exploiting their relation to the corresponding moduli spaces of 
$J$-holomorphic curves in $X$.  

\begin{defn}  \rm When $k\ge 1$,
 $\Mm_k = \Mm(A-kF,\Jj)$ is  the universal moduli space of all 
unparametrized $J$-holomorphic curves in class $A-kF$.  Thus its elements
are  equivalence classes $[h,J]$ of pairs $(h,J)$, where $J\in \Jj = \Jj^\la$, $h$
is a $J$-holomorphic map $S^2\to X$ in class $A - kF$, and where $(h,J)
\equiv (h\circ\ga, J)$ when $\ga:S^2\to S^2$ is a holomorphic
reparametrization of $S^2$.   Similarly, we write $\Mm_0=\Mm(A, x_0,\Jj)$
for the universal moduli  space of all unparametrized $J$-holomorphic curves
in class $A$ that go through a fixed point $x_0\in X$.  Thus its elements are
equivalence classes of triples $[h,z,J]$ with $z\in S^2$, $(h,J)$  as before, 
$h(z) = x_0$ and where $(h,z,J)
\sim (h\circ\ga, \ga^{-1}(z), J)$ when $\ga:S^2\to S^2$ is a holomorphic
reparametrization of $S^2$. \end{defn}

The next lemma restates part (i) of Theorem~\ref{main}.  The proof uses
standard Fredholm theory for $J$-holomorphic curves and is given  in
\S 4.1.   The only noteworthy point is that when $k > 0$ the almost
complex structures in $\Jj_k$  are not regular.  In fact, the index of the
relevant Fredholm operator is $-(4k-2)$.  However, because we are in
$4$-dimensions the Fredholm operator has no kernel, which is the basic
reason why the space of $J$ for which it has a solution is a submanifold of
codimension $4k-2$.

\begin{lemma}\label{mfld}  For all $k\ge 0$, the projection 
$$
\pi_k: \Mm_k\to \Jj_k:  \quad [h,J]\mapsto J
$$
is a  diffeomorphism of the Fr\'echet manifold $\Mm_k$ onto the
submanifold
 $\Jj_k$ of $\Jj$.  This submanifold is an open subset of $\Jj$ when $k = 0$
and has  codimension $4k-2$ otherwise. 
\end{lemma}

Our tool for understanding the stratification of $\Jj$ by the $\Jj_k$ is the
compactification $\oMm(A-kF, \Jj)$ of $\Mm(A-kF, \Jj)$ that is formed by
$J$-holomorphic stable maps. For the convenience of the reader we  recall
the definition of stable maps with $p$ marked points.  We always assume
the domain  $\Si$ to have genus $0$.  Therefore it is a connected union
$\cup_{i = 0}^m\Si_i$ of Riemann surfaces  each of which has a
given identification with the standard sphere $(S^2,j_0)$.  (Note that we
consider $\Si$ to be a topological space: the labelling of its components is a
convenience and not part of the data.) The intersection pattern of the
components can be described by a tree graph  with $m+1$
vertices, one for each  component of
$\Si$, that are connected by an edge if and only if the
corresponding components intersect.   No more than two components meet
at any point. Also,   there are $p$
 marked points $z_1,\dots, z_p$ placed anywhere on $\Si$ except at an
intersection point of two components. (Such pairs $(\Si, z_1,\dots,
z_p) = (\Si,z)$ are called semi-stable curves.) 

Now consider  a triple $(\Si, h, z)$ where  $h:\Si\to
X$ is such that $h_*([\Si]) = B$ and where the following {\it stability
condition} is satisfied:
\begin{quote}{\it  the restriction $h_i$ of the map $h$ to
$\Si_i$  is nonconstant unless $\Si_i$
contains at least  $3$ special points.}
\end{quote} 
(By definition, special points are either points of intersection  with other
components or  marked points.) 
A {\em stable map} $\si = [\Si, h,z]$ in class
$B\in H_2(X,\Z)$ is an equivalence class of such triples, where
 $(\Si, h, z')\equiv(\Si, h\circ \ga,z)$ if there is an element $\ga$ of the
group $\Aut(\Si)$ of all  holomorphic self-maps of $\Si$   such that $
\ga(z_i) = z_i'$ for all $i$.  For example, if $\Si$ has only one component and
there are no marked points, then $(\Si, h) \equiv (\Si, h\circ \ga)$ for all
$\ga\in \Aut(S^2) = \PSL(2,\C)$. Thus stable maps are {\em unparametrized}. 
We may think of the triple $(\Si, h, z)$ as a parametrized stable map. Almost
always we will only consider  stable maps that are $J$-holomorphic for
some $J$.  If necessary, we will include $J$ in the notation, writing elements
as $\si = [\Si,  h, z, J]$, but often $J$ will be understood.

Note that some stable maps $\si= [\Si, h,z,J]$ have a nontrivial
reparametrization group $\Ga_\si$.  
Given a representative $(\Si, h,z,J)$ of $\si$, this group may be defined as 
$$
\Ga_\si = \{\ga\in\Aut(\Si) : h\circ \ga = h, \ga(z_i) = z_i, 1\le i\le p\}.
$$
It is finite because of the stability condition.  The points where
this reparametrization group $\Ga_\si$ 
is nontrivial are singular or orbifold points of the moduli space.
Here is an example where it is 
nontrivial.

\begin{example}\label{smap}\rm Let $\Si$ have three components, with
$\Si_2$ and $\Si_3$ both intersecting $\Si_1$ and let $z_1$ be a marked
point on $\Si_1$.  Then we can allow 
 $h_1$ to be constant without violating stability.  If in addition $h_2,
h_3$ have the same image curve, there is an automorphism
that interchanges $\Si_2$ and $\Si_3$.  Since nearby stable maps do
not have this extra symmetry,  $[\Si, h, z_1]$ is a 
singular point in its moduli space.   However,
because  marked points are labelled, there is no such automorphism if we
put one marked point $z_2$ on $\Si_2$ and another $z_3$ at the
corresponding point  on $\Si_3$, i.e. so that $h_2(z_2) = h_3(z_3)$.  One can
also destroy this automorphism by adding just one marked point $z_0$ to
$[\Si,  h, z_1]$ anywhere on $\Si_2$ or $\Si_3$.   \end{example}

\begin{defn} \rm For $k\ge 0$ we define $\oMm(A-kF, J)$  to be the space
of all 
$J$-holomorphic stable maps $\si = [\Si, h, J]$ in class $A-kF$.  Further,
given any subset $\Kk$ of $\Jj$ we write 
$$
\oMm(A-kF, \Kk) = \cup_{J\in \Kk}\,\oMm(A-kF, J).
$$
It follows from the proof of Lemma~\ref{basic} that the domain $\Si =
\cup_{i=0}^p \Si_i$ of   $\si\in \oMm(A-kF, \Jj)$  contains a unique
component that is mapped to a curve  in some class
$A-mF$, where $m \ge k$.  We call this component the {\it stem} of $\Si$ and
label it $\Si_0$. Thus  $\oMm(A-kF, \Jj_m)$ is the moduli space of all curves
whose stem lies in class $A-mF$. Note that $\Si - \Si_0$ has a finite number
of connected components called {\it branches}.  If $h_0$ is parametrized as a
section, a branch $B_w$ that is attached to $\Si_0 = S^2$ at the point $w$ is
mapped into the fiber $\pi_J^{-1}(w)$.   In particular, distinct branches are
mapped to distinct fibers.   \end{defn}

 The moduli spaces $\oMm(A-kF, J)$ and $\oMm(A-kF, \Jj)$
have natural
stratifications, in which each stratum  is defined by fixing the topological
type of the pair $(\Si,z)$ and the homology classes
$[h_*(\Si_i)]$ of the components.  Observe that
the  class $A - mF$ of the stem is fixed  on each 
stratum $\Ss$ in $\oMm(A-kF, \Jj)$.  Hence there is a projection
$$
\Ss\to \Jj_m,
$$
whose fiber at $J\in \Jj_m$ is some stratum of $\oMm(A-kF, J)$.  
Usually, in order to have a moduli space with a nice structure one needs to
consider perturbed $J$-holomorphic curves.  But, because we are working
with genus $0$ curves in dimension $4$, the work of
Hofer--Lizan--Sikorav~[HLS] shows that  all $J$-holomorphic
curves are essentially regular.  In particular, all curves representing
some multiple $mF$ of the fiber class are regular.
Therefore each stratum of $\oMm(A-kF, J)$  is a
(finite-dimensional)  manifold. The following result is an immediate
consequence of Lemma~\ref{mfld}.

\begin{lemma}\label{strat}  Each stratum $\Ss$ of 
 $\oMm(A-kF, \Jj)$ is a manifold and the
projection $\Ss\to \Jj_m$ is a  locally trivial fibration. \end{lemma}

\begin{defn}\label{def:om}  \rm When $k \ge 1$, we set
 $\oMm_k = \oMm(A-kF,\Jj)$.  Further, $\oMm_0 = \oMm(A,
x_0, \Jj)$ is the space of all  stable maps $[\Si,  h, z, J]$ where $[\Si, h,z]$ is a 
$J$-holomorphic stable map in class $A$ with one marked
point $z$ such that $h(z) = x_0$.
\end{defn}

In the next section we show how to fit the strata of $\oMm_k$ together
 by gluing  to
form an orbifold structure on $\oMm_k$ itself.

\subsection{Gluing}

In this section we describe the structure of a neighborhood 
$\Nn(\si)\subset  \oMm_k$ of a single point $\si
\in \oMm(A-kF,  \Jj_m)$.
  Suppose that $\si = [\Si, h, J]$, and order
the components  $\Si_i$ of $\Si$ so that $\Si_0$ is the stem and so that the
union $\cup_{i\le \ell}\Si_i$ is connected for all $\ell$.  Then each $\Si_i, i>
0$ is attached to a unique component $\Si_{j_i}, {j_i} < i$ by identifying some
point $w_i\in \Si_i$ with a point $z_i\in \Si_{j_i}$.  At each such
intersection point consider the ``gluing parameter"
$$
a_i\in T_{w_i}\Si_i\otimes_{\C} T_{z_i}\Si_{j_i}.
$$
The basic process of gluing allows one to resolve the singularity of $\Si $ at
the node $w_i = z_i$ by replacing the component $\Si_i$ by a disc attached to
$\Si_{j_i}$ and suitably altering the map $h$.  As we now explain, there is a
$2$-dimensional family of ways of doing this that is parametrized by (small)
$a_i$.  

\begin{prop}\label{nbhdpt} Each $\si
\in \oMm(A-kF,  \Jj_m)$ has a neighborhood $\Nn(\si)$ 
in $\oMm_k$ that is
a product $\Uu_\Ss(\si)\times (\Nn(V_\si)/\Ga_\si)$, where 
$\Uu_\Ss(\si)\subset  \oMm(A-kF,  \Jj_m)$  is a small
 neighborhood of $\si$ in its stratum $\Ss$ and where $\Nn(V_\si)$
is a small $\Ga_\si$-invariant neighborhood of $0$ in the space of gluing
parameters $$
V_\si= \bigoplus_{i>0}\, T_{w_i}\Si_i\otimes_{\C} T_{z_i}\Si_{j_i}.
$$
\end{prop}
\proof{}  The proof is an adaptation of standard arguments in the
theory of stable maps.   The only new point is that 
the stem components are not regular so that when one
does any gluing that involves this component one has to allow $J$ to vary in a
normal slice $\Kk_J$ to the submanifold  $\Jj_m$ at $J$.  This analytic detail is
explained in \S4.2.
What we will do here is describe the topological aspect of the proof.

First of all, let us describe the process of gluing.
Given $a\in V_\si$, the idea is first to
construct an approximately $J$-holomorphic  stable map $(\Si_a, h_a, J)$ on a
glued domain $\Si_a$ and then to perturb $h_a$ and $J$ using a Newton
process to a $J_a$-holomorphic map 
 $h_a: \Si_a\to X$ in $\oMm(A-kF, \Kk_J)$. We will describe the first step
in some detail here since it will be used in \S 3.  The analytic arguments 
needed for the second step are postponed to \S 4.

The glued domain
$\Si_a$ is constructed as follows.  For each $i$ such that $a_i \ne 0$, cut  out 
 a small open disc ${\INT} D_{w_i}(r_i)$  in $\Si_i$ centered at  $w_i$ and a
similar disc ${\INT} D_{z_i}(r_i)$ in  $B_{j_i}$ where
 $r_i^2=\|a_i\|$, and then  glue the boundaries of these discs together
with a twist prescribed by  the argument of $a_i$.  The Riemann surface
$\Si_a$ is the result of performing this
operation for each $i$ with $a_i\ne 0$.  
(When $a_i = 0$ one simply leaves the component $\Si_i$ alone.)

To be more precise, consider  gluing   $z\in \Si_0$ to $w\in
\Si_1$.   Take a K\"ahler metric on $\Si_0$ that is flat near  $z$ and
 identify the disc $D_z(r)$ isometrically with the disc of radius $r$ in the
tangent space $T_z = T_z(\Si_0)$ via the exponential map.  Take a similar
metric on $(\Si_1,w)$. Then the gluing $\p
D_z(r)\to \p D_w(r)$ may be considered as the restriction of the map 
$$
\Psi_a: T_z - \{0\} \;\longrightarrow\; T_w - \{0\} 
$$
 that is defined for $x\in T_z$ by the requirement: $$
x\otimes \Psi_a(x) =  a,\quad x\in T_z.
$$
Thus, with respect to  chosen identifications of $T_z$ and $T_w$ with $\C$,
$\Psi_a$ is given by the formula: $x\mapsto  a/x$ and so takes the circle of
radius $r = \sqrt{\|a\|}$ into itself.   This describes the glued domain $\Si_a$
as a point set. It remains to put a metric on $\Si_a$ in order to
 make it a Riemann surface.  By hypothesis the original metrics on $\Si_0,
\Si_1$ are flat near $z$ and $w$ and so may be identified with  the flat
metric  $|dx|^2$ on $\C$.  Since 
$$
\Psi_a^*(|dx|^2) = \left|\frac a{x} \right|^2|dx|^2,
$$
$\Psi_a(|dx|^2) = |dx|^2$ on the circle $|x| = r$.  Hence,  we may choose a
function
$
\chi_r:(0,\infty)\to (0,\infty)$ so that the metric $\chi_r(|x|)|dx|^2$ is
invariant by $\Psi_a$ and so that $\chi_r(s) = 1$ when $s > (1+\eps)r$, and
then patch together the given metrics on $\Si_0 - D_z(2r)$ and $\Si_1 -
D_z(2r)$ via  $\chi_r(|x|)|dx|^2$.

In \S 3 we  need to understand what happens as $a$
rotates around the origin. It is not hard to check  that if we write $a_\th =
e^{i\th} a_z\otimes a_w$, where $a_z\in T_z, a_w\in T_w$ are fixed, and if
$\Psi_{a_\th}$ identifies the point $p_z$ on $\p D_z(r)$ with $p_w$ on $\p
D_w(r)$ then $$
\Psi_{a_\th} (e^{i\th} p_z) = p_w.
$$

The next step is to define the approximately holomorphic map (or
pre-gluing) $h_a: \Si_a\to X$ for sufficiently small $\|a\|$.  The map
$h_a$ equals $h$ away from the discs $D_{z_i}(r_i), D_{w_i}(r_i)$, and
elsewhere is defined by using cut-off functions that  depend only on
$\|a\|$.  To describe the deformation of $h_a$ to a holomorphic map one 
needs to use analytical arguments.  Hence further details are postponed until
\S 4.  

We are now in a position to describe a neighborhood of $\si$. 
It is convenient to think of $V_\si$ as the direct sum $V_\si' \oplus V_\si''$
where $V_\si'$ consists of  the summands $T_{w_i}\Si_i\otimes_{\C}
T_{z_i}\Si_{j_i}$ with $j_i = 0$ and $V_\si''$ of the rest.  Note that the
obvious action of $\Ga_\si$ on $V_\si$ preserves this splitting. 
(It is tempting to think that the induced action on $V_\si'$ is trivial since the
elements of $\Ga_\si$ act trivially on the stem.  However, this need not
be so since they may rotate branch components that are attached to the stem.)
If we glue at points parametrized by $a''\in V_\si''$ then the corresponding
curves lie in some branch and are regular.  Hence the result of gluing is a
$J$-holomorphic curve (i.e. there is no need to perturb $J$).   Further, because
the gluing map $\TGg$ is $\Ga_\si$-equivariant, there is a
neighborhood    of $\si$ in $ \oMm(A-kF,  \Jj_m)$ of the form
$$
\Uu''(\si) = \Uu_\Ss(\si)\times (\Nn(V_\si'')/\Ga_\si),
$$
 where $\Nn(V)$ denotes
a neighborhood of $0$ in the vector space $V$.

When we glue with
elements from $V_\si'$, the homology class of the stem changes and so
the result cannot be $J$-holomorphic since $J\in \Jj_m$.  We show in
Proposition~\ref{exun}
that if $\Kk_J$ is a normal slice to  the submanifold $\Jj_m$ at $J$ then 
for sufficiently small $a\in V_\si'$ the
approximately holomorphic map $h_a:\Si_a\to X$ deforms to a
unique $J_a$ holomorphic map $\TGg(h_\si, a)$ with $J_a\in \Kk_J$. 
Therefore, for each element $\si'' =
[\Si, h'',J''] \in \Uu''(\si)$  there is a homeomorphism from some
neighborhood  $\Nn(V_{\si''}')$ onto a neighborhood of $\si''$ in 
$\oMm(A-kF, \Kk_{J''})$.  Moreover, if $\Uu''(\si)$ is sufficiently small,  the
spaces $V_{\si''}$ can all be identified with $V_\si$ and it follows from
the proof of Proposition~\ref{exun} that the neighborhoods $\Nn(V_{\si''})$
can be taken to have uniform size and so may all be identified.  
Hence  the
neighborhood $\Nn(\si)$ projects to $\Uu''(\si)$ with fiber at $\si''$ equal
to $\Nn(V_\si')/\Ga_{\si''}$.  
 In general, the groups $\Ga_{\si''}$ are
subgroups of $\Ga_\si$ that vary with $\si''$: in fact they equal the stabilizer
of the corresponding gluing parameter $a''\in V(\si'')$.   However, since
elements of $\Uu_\Ss(\si)$ lie in the same stratum they have isomorphic
isotropy groups.   It is now easy to
check that the composite map
$$
\Nn(\si)\to \Uu''(\si)\to \Uu_\Ss(\si)
$$
 has fiber $\Nn(V_\si)/\Ga_\si$
as claimed.
\QED

\subsection{Moduli spaces and the stratification of $\Jj$}

Since each stable $J$-curve in class $A - kF$ has exactly one component
in some class $A- mF$ with $m \ge k$, the
 projection $\pi_k: \oMm(A-kF, \Jj)\to \Jj$ has image $\oJj_k$.
Consider the inverse image
$$
\oMm(A-kF, \Jj_m) = \pi_k^{-1}(\Jj_m).
$$
The next result shows that we can get a handle on the structure of
$\oJj_k$ by looking at the spaces $\oMm(A-kF, \Jj_m)$.

\begin{prop}\label{prop:fibk} When $k > 0$ the projection
$$
\pi_k: \oMm(A-kF,  \Jj_m)\to \Jj_m
$$
is a locally trivial fibration whose fiber $\Ff_J(m-k)$ at $J$
is  the space of all stable $J$-curves $[\Si,h]$ in class $A-kF$ that have as
one component the unique $J$-holomorphic curve $\Delta_J$ 
in class $A-mF$.  In
particular, $\Ff_J(m-k)$ is  a stratified space with strata that are manifolds
of (real) dimension $\le 4(m-k)$. Its
 diffeomorphism type  depends only on $k-m$.
\end{prop}
\proof{}  Let us look at the structure of $\Ff_J(m-k) =
\pi_k^{-1}(J)$. The stem of each element $[\Si, h,J]\in \Ff_J(m-k)$  is mapped
to the unique $J$-curve $\Delta_J$ in class $A-mF$.  Fix this component
further by supposing that it is parametrized as a section of the fibration
$\pi_J:X\to S^2$ (where $\pi_J$ is as in Lemma~\ref{basic}.)  
We may divide the fiber $\Ff_J(m-k)$ into disjoint sets $\Zz_{\Dd,J}$ each
parametrized by  a fixed  decomposition $\Dd$
of $m-k$ into  a sum $d_1 + \dots + d_p$ of unordered positive numbers.  
The elements of $\Zz_{\Dd,J}$ are those with $p$ branches 
$B_{w_1},\dots, B_{w_p}$ where $h_*[B_{w_i}] = d_i [F]$.  Thus $\Zz_{\Dd,J}$
maps onto the configuration space of $p$ distinct (unordered) points in $S^2$
labelled by the positive integers  $d_1,\dots, d_p$ with sum $m-k$. 
Moreover this map is a fibration with fiber equal to the product
$$
\prod_{i=1}^d \oMm_{0,1}(S^2, q, d_i)
$$
where $\oMm_{0,1}(S^2, q, d)$ is the space of $J$-holomorphic
stable maps into $ S^2$ of degree $d$ and with one marked point $z$ such
that $h(z) = q$. (This point $q$ is where the branch is attached to $\De_J$.) 
According to the general theory, $\oMm_{0,1}(S^2, q, d)$ is an orbifold of real
dimension $4(d-1)$.  It follows easily that $\Zz_{\Dd,J}$ is an orbifold of real
dimension $4(m-k)-2p$.  It remains to understand how the different sets
$\Zz_{\Dd,J}$ fit together, i.e. what happens when two or more of the points
$w_i$ come together.  This may be
described by suitable gluing parameters as in Proposition~\ref{nbhdpt}.  
The result follows.  (For more details see any reference on stable maps,
eg [FO], [LiT], [LiuT1].  An example is worked out in \S3.2.4 below.) \QED

\NI
{\bf Note}  For an analogous statement  when $k=0$ 
see Proposition~\ref{prop:fib0}.
\MS

Our next aim is to describe the structure of a neighborhood of
$\oMm(A-kF,  \Jj_m)$ in $\oMm_k=\oMm(A-kF, \Jj)$.
We will write $\Zz_J$ for the fiber $\Ff_J(m-k)$ of $\pi_k$ that was
considered above and set
$$
\Zz= \bigcup_{J\in \Jj_m}\, \Zz_J,\qquad \Zz_\Dd = \bigcup_{J\in \Jj_m}\, 
\Zz_{\Dd,J}.
$$
(The letter $\Zz$ is used here because $\Zz$ is the ``zero-section" of the space
of gluing parameters $\Vv$ constructed below.)
 Consider an element $\si = [\Si, h, J] $ that lies in a
substratum $\Zz_\Ss$ of $ \Zz_\Dd$ where $\Dd = d_1+\dots + d_p$.  Then
$\Si$ has $p$ branches $B_1,\dots, B_p$  that are attached at the distinct
points $w_1,\dots, w_p\in \Si_0$.  Let $z_i$ be the point in $ B_i$ that is
identified with $w_i\in \Si_0$ and define
$$
V_{\si} = \bigoplus_{i=1}^p \; T_{z_i}B_i\otimes_{\C} T_{w_i}\Si_0.
$$
As explained in Proposition~\ref{nbhdpt}, the gluing parameters $a\in V_\si$
(when quotiented out by $\Ga_\si$) parametrize a normal slice to
$\Zz_\Dd$ at $\si$.  (Note that previously $V_\si$ was called $V_\si''$.)
We now want to show how to fit these vector spaces  together to 
form the fibers of an orbibundle\footnote
{
A rank $k$ orbibundle $\pi: E\to Y$ over an orbifold $Y$ has the following
local structure.   Suppose that $\si\in Y$ has local chart
$U\subset \TU/\Ga_\si$ where the uniformizer $\TU$ is a subset of
$\R^n$.  Then $\pi^{-1}(U)$ has the form $\TU\times \R^k/\Ga_\si$ where
the action of $\Ga_\si$ on $\R^n\times \R^k$ lifts that on $\R^n$ and is
linear on $\R^k$. There is an obvious compatibility condition between
charts: see~[FO],\S2.} 
over $\Zz_\Dd$.  Here we must  incorporate twisting that arises from the
 fact that gluing takes place on the space of {\it parametrized }
stable maps. Since this is an important point, we dwell on it at some length.
For the sake of clarity, we will in the next few paragraphs denote parametrized
stable maps by $\tsi = (\Si,h)$ and the usual (unparametrized) maps by
$\si = [\Si, h]$.  Further, $\Ga_\tsi$ denotes the corresponding realization of the
group $\Ga_\si$ as a subgroup of $\Aut(\Si)$.

Recall that $X$ is identified with $S^2\times S^2$ in such a way  that the
fibration $\pi_J: X\to S^2$ whose fibers are the $J$-holomorphic $F$-curves is
simply given by projection onto the first factor.  Hence each such fiber has a
given identification with $S^2$.  Further, we  assume that the stem
$h_{\si,0}:\Si_{\si,0}\to \Delta_J$ is parametrized as a section $z\mapsto (z,
\rho(z))$.  Hence  we only have to choose parametrizations of each branch.  
Since each branch component has at least one special point, its
automorphism group is either trivial or has the homotopy type of $S^1$.
Let ${\rm Aut\,}'(\Si)$ be the subgroup of ${\rm Aut\,}(\Si)$ consisting of
automorphisms that are  the identity on the stem. Then the identity
component of  ${\rm Aut\,}'(\Si)$  is homotopy equivalent to  a torus
$T^{k(\Ss)}$.  (Here $\Ss$ is the label for the stratum containing $\si$.)
   Let $g$  be
 a $\Ga_\tsi$-invariant metric on the domain $\Si$ that is  also
invariant under some action of the torus $T^{k(\Ss)}$.  

\begin{defn}\label{def:gp}\rm  The group $\Aut^K(\Si)$ is defined to be the
subgroup of  the isometry group of $(\Si, g)$ generated by $\Ga_\tsi$ and
$T^{k(\Ss)}$.   Note that  $\Ga_\tsi$  is the semidirect
product of a subgroup $\Ga_\tsi'$ of $T^{k(\Ss)}$ with a subgroup $\Ga_\tsi''$
that permutes the components of each branch.  Further  $\Aut^K(\Si)$ is a
deformation retract of the subgroup $p^{-1}(\Ga_\tsi'')$ of $\Aut(\Si)$, where
we consider $\Ga_\tsi''$ as a subgroup of $\pi_0(\Aut(\Si))$ and
$$
p:\Aut(\Si)\to \pi_0(\Aut(\Si))
$$
is the projection.  For a further discussion, see \S4.2.4.
\end{defn}

Let us first consider a fixed $J\in \Jj_k$.
It follows from the above discussion that
 on each stratum $\Zz_{\Ss,J}$ there is a principal bundle
$$
\Zz_{\Ss,J}^{para} \to \Zz_{\Ss,J}
$$
with fiber $\Aut^K(\Si)$ such that the elements
of $\Zz_{\Ss,J}^{para}$ are parametri\-zed stable maps $\tsi = (\Si, h)$. 
 Since the  space $V_\tsi$ of gluing parameters at $\tsi$ is
made from tangent spaces to $\Si$ there is a well defined bundle
$$
\Vv_{\Ss,J}^{para} \to \Zz_{\Ss,J}^{para}
$$
with fiber $V_\tsi$.  Further the action of the reparametrization
group $\Aut^K(\Si)$ lifts to $\Vv_{\Ss,J}^{para}$  and we define $\Vv_{\Ss,J}$ to
be the quotient  $\Vv_{\Ss,J}^{para}/\Aut^K(\Si)$.  Thus there is a commutative
diagram $$
\begin{array}{ccc} \Vv_{\Ss,J}^{para} & \to & \Vv_{\Ss,J}\\
\downarrow &  & \downarrow\\
\Zz_{\Ss,J}^{para} & \to & \Zz_{\Ss,J}.
\end{array}
$$
where the right hand vertical map is an orbibundle with fiber
$V_\tsi/\Ga_\tsi$.

Now consider the space
$
\Zz_{\Dd,J} = \cup_{\Ss\subset \Dd} \Zz_{\Ss,J}.
$
 The local topological structure of
$\Zz_{\Dd,J}$ is given by
gluing parameters as in  Proposition~\ref{nbhdpt}.  Observe that
every $J$ is regular for the branch components so that
the necessary gluing operations can be performed keeping $J$ fixed.
The spaces $\Zz_{\Dd,J}^{para},\Vv_{\Dd,J}^{para}$  are defined similarly  and
clearly there is a  vector bundle $\Vv_{\Dd,J}^{para}\to \Zz_{\Dd,J}^{para}$,

We want to see that the union 
$$
\Vv_{\Dd,J} = \bigcup _{\Ss\in \Dd} \Vv_{\Ss,J}.
 $$
has the structure of an orbibundle over
$\Zz_{\Dd,J}$. The point here is that the groups $\Aut^K(\Si)$ change 
dimension as $\tsi$
moves from stratum to stratum.  Hence we need to see that the local gluing
construction that fits   the different strata in $\Vv_{\Dd,J}$   together is
compatible with the group actions.   We show in \S4.2.4 below
that the gluing map $\TGg$ can be defined at the point $\tsi$ to be
$\Aut^K(\Si)$-invariant, i.e. so that 
$$ 
\TGg(h_\si, a) = \TGg(h_\si\circ
\theta^{-1}, \theta\cdot a), 
$$
where $\TGg(h_\si,a)$ is the result of gluing the map $h_\si$ with
parameters $a$.
In the situation considered here, we are dividing the set of gluing parameters
at $\tsi$ into two, and will write $a = (a_b, a_s)$ where $a_b$ are the gluing
parameters at intersections of branch components and $a_s$ are those
involving the stem component.  As $h_\si$ moves within $\Zz_{\Dd,J}^{para}$
we glue along $a_b$, considering  $a_s$ to be part of the fiber
$V_\si$.   Moreover, 
if $\tsi'  = (\Si_{a_b}, \TGg(h_\si,
a_b))$, Lemma~\ref{le:repr} (ii)  shows that 
 $\TGg$ can be constructed to
be compatible with the actions of the groups
$\Aut^K(\Si)$ and $\Aut^K(\Si')$ on the fibers $V_\si$ and
$V_{\si'}$ of $\Vv_{\Dd,J}^{para}$ at $\tsi,\tsi'$.  
It follows without difficulty that the quotient 
$$
\Vv_{\Dd,J}\to \Zz_{\Dd,J}
$$
is an orbibundle. 

Finally, one forms spaces
$$
\Vv_J = \bigcup_\Dd\;\Vv_{\Dd,J},\qquad \Vv = \bigcup_{J\in \Jj_k}\;\Vv_J 
$$ 
whose local structure is also described by appropriate gluing parameters 
as above.
Forgetting the gluing parameters gives projections
$$
\Vv_J\to \Zz_J = \Ff_J(m-k);\qquad \Vv\to \Zz = \Ff(m-k),
$$
and $\Zz_J, \Zz$ embed in
$\Vv_J$ and $\Vv$  as the ``zero sections".  The map $\Vv_J\to
\Zz_J$ preserves the stratifications of both spaces.  However it is no longer an
orbibundle  since the dimension of the fiber $V_\si$ depends on
$\Dd$.  In fact, the way that the different sets $\Vv_{\Dd,J}$ are fitted
together is best thought of as a kind of plumbing: see \S3.2.4.

\begin{example}\label{ex:lens}\rm Everything is greatly simplified when $m
- k = 1$.  Here there is only one decomposition $\Dd$ and the space
$\Zz_{\Dd,J}$ consists of just one stratum  diffeomorphic to  $S^2$. 
Moreover the bundle $\Zz_{\Dd,J}^{para}\to \Zz_{\Dd,J}$ has a section
with the following description.  Choose $J\in \Jj_{k+1}$ so that $\pi_J:X\to S^2$
is the standard projection onto the first factor and so  that the graph $h_0$ of
the map $\rho_k:S^2\to S^2$ of degree $-(k+1)$ is $J$-holomorphic.  Let $\Si_0,
\Si_1$ be two copies of $S^2$ and for each $w\in S^2$ define $(\Si_w, h_w)\in
\Zz_{\Dd,J}^{para}$ by   
\begin{eqnarray*} 
\Si_w & = & \Si_0\, \cup_{w = \rho(w)}\,
\Si_1,\\ 
h_w|_{\Si_0} = h_0,& &  h_w|_{\Si_1}: z\mapsto (w,z).
\end{eqnarray*}
Hence in this case $\Vv_J$ is a complex line bundle over $\Zz_J = S^2$. To
calculate its Chern class, observe that $\Vv_J$ can be identified with the space
$$
\bigcup_{w\in S^2} T_{\rho(w)}\Si_1\otimes T_w(\Si_0) =
(TS^2)^{-k-1}\otimes TS^2,
$$
and so has Chern class $-2k$.
\end{example}

The following result is proved in \S4.

\begin{prop}\label{glue} There is a neighborhood $\Nn_\Vv(\Zz)$ of  $\Zz$
in $\Vv$ and a gluing map 
$$
\Gg: \Nn_\Vv(\Zz) \longrightarrow \oMm(A-kF, \Jj)
$$ 
that maps $\Nn_\Vv(\Zz)$  homeomorphically onto a neighborhood of 
$\oMm(A-kF,\Jj_m)$ in $\oMm(A-kF, \Jj)$.
\end{prop}

It follows from the construction of $\Gg:\Nn_\Vv(\Zz)\to \oMm(A-kF,\Jj_m)$ 
outlined in Proposition~\ref{nbhdpt}  that the stem of the glued map $\Gg(\si,
a)$ lies in the class $A - (m - \sum_i n_i)F$ where the indices $i$ label the
branches $B_i$ of $\si$ and $n_i$ is defined as follows.  
 If $a_i = 0$ then $n_i = 0$.  Otherwise,  if $\Si_{j_i}$ is
the component  of $B_i$ that meets $\Si_0$ then 
 $n_i$ is the multiplicity of $h_{j_i}$, that is $[h(\Si_{j_i})] = n_i F$.  
Let $\Nn_p$ denote the set of all elements $(\si,a)\in \Nn_\Vv(\Zz)$ such
that the stem of the glued map lies in class $A-pF$.  In other words,
$$
\Nn_p = (\pi_k\circ\Gg)^{-1} \Jj_p.
$$
Clearly, $\Nn_p$ is a union of strata in the stratified space $\Nn_\Vv(\Zz)$.
Further, when $k > 0$ the map $\pi_k: \Gg(\Nn_p)\to \Jj_p$ is a
fibration with fiber $\Ff(p-k)$.  The next proposition follows immediately
from Proposition~\ref{prop:fibk}.

\begin{prop}\label{link1}  The link $\Ll_{m,k}$ is the finite-dimensional
stratified space obtained from the link of $\Zz_J$ in $\Vv_J$
 by collapsing the fibers  of the projections $\Vv_J\cap \Nn_p
\to \Jj_p$ to single points.
\end{prop}
       
\NI
{\bf Proof of Proposition~\ref{pr:lens}}  

We have to show that the link $\Ll_{k+1,k}$ is the lens space $L(2k,1)$.
We saw in
 Example~\ref{ex:lens} that $\Vv_J$ is a line bundle with Chern class $-2k$.
In this case there is only  one nontrivial
stratum in $\Nn_\Vv(\Zz)$, namely $\Nn_k$, which is the complement of the
zero section.  Moreover, the map $\pi_k\circ\Gg$ is clearly injective. 
Hence by the above lemma $\Ll_{k+1,k}$ is simply the unit sphere bundle
of $\Vv_J$ and so is a lens space as claimed. \QED

\section{The link $\Ll_{2,0}$ of $\Jj_2$ in $\oJj_0$}

In this  section we illustrate Proposition~\ref{link1} by calculating the link
$\Ll_{2,0}$.  We know from Lemma~\ref{mfld} that $\Ll_{2,0} = S^5$. 
 The general theory of \S2 implies
that $\Ll_{2,0}$ can be obtained from the link $\Ll_{\Zz}$ of the zero section
$\Zz_J$ in the stratified space $\Vv_J$ of gluing data by collapsing certain
strata.  When looked at from this point of view, the
$S^5$ appears in quite a complicated way that was described in 
Theorem~\ref{LINK}.
We  begin here by explaining the
plumbing construction, and then discuss how this relates to $\Ll_\Zz$.

\subsection{Some topology}
 
Recall that $S(L_P)\to \Pp (\Oo(k)\oplus \Oo(m))$ is the unit circle
bundle of the canonical line bundle $L_P$ over the projectivization $\Pp
(\Oo(k)\oplus \Oo(m))$.

 \begin{lemma}\label{le:str}
The bundle $S(L_P) \to \Pp(\Oo(-1)\oplus
\C)$ can be identified with the pullback of the canonical circle bundle
$S(L_{can})\to \C P^2$ over the blowdown map
$\C P^2\# \bcp\to \C P^2$.
\end{lemma}
\proof{} It is well known that $\Pp(\Oo(-1)\oplus
\C) $ can be identified with $\C P^2\# \bcp\to \C P^2$.  Indeed the section 
$S_- = \Pp(\{0\}\oplus\C)$ has self-intersection $-1$, while $S_+ =
\Pp(\Oo(-1)\oplus \{0\})$  has self-intersection $1$. 
Further, the circle bundle $L_P$ is trivial over $S_-$ and has Euler class 
$-1$ over $S_+$ and over the fiber class. The result follows.\QED

The space we are interested in is formed by plumbing a rank $2$ bundle
$E\to S^2$ to a  line bundle $L\to Y$, where $\dim (Y) = 4$.  This plumbing
$E\pl L$ is the space
obtained from the unit disc bundles $D(E)\to S^2$ and $D(L)\to Y$ by 
identifying the inverse images of discs $D^2, D^4$ on the two bases
in the obvious way:
 the disc fibers of $D(E)\to S^2$ are identified with
with flat sections of $D(L)$ over $D^4$ and flat sections of $D(E)$
over $D^2$ are identified
with  fibers of $D(L)$.  There is a corresponding plumbing 
$S(E)\pl S(L)$ of the
two sphere bundles $S^3\to S(E)\to S^2$ and $S^1\to S(L)\to Y$, obtained by
cutting out the inverse images of open discs in the two bases and 
appropriately gluing the boundaries. 
The resulting space  $S(E)\pl S(L)$ is the link of the core $S^2\cup Y$ in 
the plumbed bundle $E\pl L$.

\begin{lemma}  Let $L_{can}\to \C P^2$ be the canonical line bundle and $E =
\Oo(k)\oplus \Oo(m)$. Then $E\pl L_{can}$ may be identified with the
blow-up of $\Oo(k+1)\oplus \Oo(m+1)$ at a point on its zero section. Hence  
$$
S(E)\pl S(L_{can}) = S((\Oo(k+1)\oplus
\Oo(m+1)).
$$
  \end{lemma}
\proof{}    First consider the structure of the blow up 
$\widetilde{\C}\,\!^3$ of
$\C^3 = \C\times \C^2$ at the origin.  The fibration $\pi: \C
\times \C^2\to
\C$ induces a fibration 
$$
\widetilde{\pi}: \widetilde{\C}\,\!^3\to \C.
$$
Clearly,  the inverse image $\widetilde{\pi} \,\!^{-1}(z)$ of each point
$z\ne 0$ is a copy of $\C^2$ while 
$\widetilde{\pi} \,\!^{-1}(0)$ is the union of the exceptional divisor together
with the set of lines in the original fiber $\pi^{-1}(0)$.  Let $\la$ be the line
 in $\C\times \C^2$ through the origin and the point $(1,a,b)$. Lift $\la$ to
the blow-up and consider its intersection with 
$$
\widetilde{\pi} \,\!^{-1}(S^1) = {\pi}^{-1}(S^1) \subset \C\times \C^2,
$$
 where 
$S^1$ is the unit circle in $\C$. This intersection consists of the points
$
(e^{it}, e^{it}a,e^{it}b)$, hence it is these circles (rather than the circles
$(e^{it}, a,b)$) that bound discs in the blowup.  
Therefore, if we think of the
blowup $\widetilde{\C}\,\!^3$ as the plumbing of the bundle
$\pi:\C\times \C^2\to \C$ with $L_{can}$,  the original trivialization of $\pi$
differs from the trivialization (or product structure) near $\pi^{-1}(0)$ that is
used to construct the plumbing.   

Now recall that 
$$
\Oo(k) = D^+\times \C \cup_\alpha D^-\times \C,
$$
where $D^+, D^-$ are $2$-discs, with $D^+$ oriented positively and $D^-$
negatively, and where the gluing map $\alpha$ is given by
$$
\alpha: \p D^+\times \C\to \p D^-\times \C:\quad (e^{it}, w)\mapsto 
(e^{it},e^{-ikt}w).
$$
It follows easily that the blowup of $D(\Oo(k+1)\oplus
\Oo(m+1))$ at a point on its zero-section is obtained by plumbing the disc
bundle $D(\Oo(k)\oplus \Oo(m))$ with
$D(L_{can})$.  This proves the first statement.  The second statement is then
immediate. \QED

We are interested in plumbing not with $L_{can}\to \C P^2$ but with
a particular singular line bundle (or orbibundle) $L_Y\to Y$.  This means
that the unit circle bundle $S(L_Y)\to Y$ is a Seifert fibration with 
a finite number of  
singular (or multiple) fibers.  In our case, there is an $S^1$ action on $S(L_Y)$
such that the fibers of  the map $S(L_Y)\to Y$ are the $S^1$-orbits.  In fact,
we can identify $S(L_Y)$ with $S^5$ in such a way that the $S^1$ action is
$$
\th\cdot (x,y,z) = (e^{i\th}x, e^{i\th}y, e^{2i\th}z),\quad x,y,z\in \C.
$$
Thus there is one singular fiber that goes through the point $(0,0,1)$.  All
other fibers $F$ are regular.  For each such $F$ there is a diffeomorphism of
$S^5$ that takes $F$ to the circle $ \ga_0=(e^{i\th}, e^{i\th}, 0)$.  Identify a
neighborhood of $\ga_0$ with $S^1\times D^4$ in such a way that
$$
S^5 = S^1\times D^4\cup D^2\times S^3,
$$
with the identity map of $S^1\times S^3$ as gluing map. Then, in these
coordinates near $\ga_0$ the fibers of $S(L_Y)$ are (diffeomorphic to) 
 the circles  
$$
\ga_x = \left\{(\th, A_{\th}(x))\in S^1\times D^4:  A_\th =
\left(\begin{array}{cc}e^{i\th} &0\\0&e^{2i\th}\end{array}  \right)\right\}.
$$
By way of contrast, the fibers of $S^5$ with the Hopf fibration have
neighborhoods fibered by the circles
$$
\ga_x' = \left\{(\th, A_{\th}'(x))\in S^1\times D^4:  A_\th' =
\left(\begin{array}{cc}e^{i\th} &0\\0&e^{i\th} \end{array} \right)\right\}.
$$

The next result shows that plumbing with $S(L)$ is a kind of twisted
blowup.

 \begin{prop}\label{str} Let $L_Y\to Y$ be the orbibundle described 
in the previous paragraph.
Then the manifold obtained by plumbing $S(\Oo(k)\oplus \Oo(m))$ with a
regular fiber of $S(L_Y)$  is diffeomorphic to $S(\Oo(k+2) \oplus
\Oo(m+1))$. \end{prop}
 \proof{}  We may think of plumbing as the result of a surgery that matches
the flat circles $S^1\times pt$ in the copy of $S^1\times S^3$ in 
$S(\Oo(k)\oplus \Oo(m))$ with the circles $\ga_x$ in the neighborhood of a
regular fiber $\ga_0$ of $S(L_Y)$.  We would get the same result if we
matched  the circles
$$
\de_x = \left\{(\th, A_{\th}''(x))\in S^1\times S^3:  A_\th'' =
\left(\begin{array}{cc}e^{-i\th} &0\\0& 1 \end{array} \right)\right\}
$$
in $S^1\times S^3\subset S(\Oo(k)\oplus \Oo(m))$ with the circles $\ga_x'$
in the standard (Hopf) $S^5$.  But  if we trivialize the
boundary of $S(\Oo(k)\oplus \Oo(m)) - D^2\times S^3$ by the circles $\de_x$
we get the same as if we trivialized the boundary of 
$S(\Oo(k+1)\oplus \Oo(m)) - D^2\times S^3$ in the usual way by flat circles.
Thus
\begin{eqnarray*}
S(\Oo(k)\oplus \Oo(m))\pl S(L_Y) & = & S(\Oo(k+1)\oplus \Oo(m))\pl
S(L_{can}) \\ &= & S(\Oo(k+2)\oplus \Oo(m+1)).
\end{eqnarray*}
There is a question of orientations here: do we have to add or subtract $1$
from $k$ to compensate for the extra twisting in $S(L_Y)$?   One can
check that it is correct to add $1$ by using the present approach to give
 an alternate proof of the previous lemma.  For, if we completely
untwisted the circles  in the neighborhood of $\ga_0$ (thereby increasing
the twisting of the other side by an additional $2$), we would be doing the
trivial surgery in which the attaching map is the identity.  Note also that
because the sum $\Oo(k)\oplus \Oo(m)$ depends only on $k+m$ we
could equally well have put the extra twist on the other factor. 
\QED

\subsection{Structure of the pair $(\Vv_J, \Zz_J)$}

Our aim  is to prove the following proposition, where
 $\Vv_J$ is the space of gluing parameters
for a fixed $J\in \Jj_2$ that describes the link of the space of
$(A-2F)$-curves in the space of (pointed) $A$-curves.

\begin{prop}\label{prop:link} The link $\Ll_{\Zz}$ of the zero section
$\Zz_J$ in the stratified space $\Vv_J$ is constructed by plumbing
$S(\Oo(-3)\oplus \Oo(-1))$ to $S(L_Y)$.  Hence
$$
\Ll_{\Zz} = S(\Oo(-1)\oplus \C).
$$
\end{prop}

 We are now not quite in the situation described in
Proposition~\ref{prop:fibk} because we are including the open stratum 
$\Jj_0$ of
$\Jj$.  This means that we have to replace the space $\oMm_k = \oMm(A-kF,
\Jj)$ by a space $\oMm_0$ of curves of class $A$ that go through the fixed
point $x_0$.  Since we are interested in working out the structure of the 
fiber of the projection $\oMm_0\to\Jj$ at a point $J\in \Jj_2$, we will choose
$x_0$ so that it does not lie on the unique $J$-holomorphic $(A-2F)$-curve
$\Delta_J$ and then define $\oMm_0$ to be the space $\oMm(A,x_0, \Jj)$ in
Definition~\ref{def:om}.  Let $\pi_0$ denote the projection  
$$
\pi_0:\oMm_0\to\Jj
$$
and set $\oMm_0(\Jj_m) = \pi_0^{-1}(\Jj_m)$ as before.  It is not hard to see
that the following analog of Proposition~\ref{prop:fibk} holds.

\begin{prop}\label{prop:fib0} (i)  Let $J\in \Jj_m$ be any 
almost complex structure such that
the unique $J$-holomor\-phic $(A-mF)$-curve $\Delta_J$ does not go through
$x_0$. Then the projection 
$$
\pi_k: \oMm_0(\Jj_m)\to \Jj_m
$$
is a locally trivial fibration near $J$, whose fiber $\Ff_J(0,m)$  
is  the space of all stable $J$-curves $[\Si,u]$ in class $A$ that have
 $\Delta_J$   as
one component and go through $x_0$.  In particular, $\Ff_J(0,m)$ is  a
stratified space whose strata are
 orbifolds of (real) dimension $\le 4m -2$.
\SS

\NI
(ii)  The singular fibers of $\pi_k: \oMm_0(\Jj_m)\to \Jj_m$ occur at points $J$
for which  $x_0\in \De_J$.  For such $J$, $\pi_k^{-1}(J)$
can be identified with the
space $\Ff_J(m)$ described in Proposition~\ref{prop:fibk}.
 \end{prop}

As before, we now construct a pair $(\Vv_J, \Zz_J)$ that describes a
neighborhood of $\oMm_0(\Jj_m)$ in $\oMm_0$.  We will concentrate on
the case $m = 2$ and will suppose that $x_0\notin \De_J$.  We further
normalize $J$ by requiring  that
 the projection $\pi_J$ along the $J$-holomorphic $F$-curves is
simply the projection onto the first factor $S^2$.   We write $q_0 = 
\pi_J(x_0)$.

\MS

\NI
{\bf \S 3.2.1  The bundle  $\Vv_{2,J}\to \Zz_{2,J}$.}\SS

Observe first that  $\Zz_J$ has
two  subsets: $\Zz_{1,J}$ consisting of all stable $A$-maps  $[\Si, z_0, h]$ 
that are the union of the $(A-2F)$-curve $\Delta_J$ with a double covering of
the fiber $F_0$ through $x_0$, and $\Zz_{2,J}$  consisting of all stable $A$-maps
$[\Si, z_0, h]$  that are the union of $\Delta_J$ with two distinct fibers.  
We will call these sets $\Zz_{i,J}$ strata.  This is accurate as far as $\Zz_{2,J}$ is
concerned, but strictly  speaking $\Zz_{1,J}$ is a union of strata.
(Recall that the strata are detemined by the topological type of the
marked domain $[\Si, z_0]$ , and the homology class of the images of its
components under $h$.)
 
Let us first consider $\Zz_{2,J}$.  Since
$h(z_0) = x_0$ always, one of the two fibers has to be $F_0$ and the other
moves.   Therefore,  the stratum $\Zz_{2,J}$ maps onto
$S^2 - \{q_0\}$.  It is convenient to compactify $\Zz_{2,J}$ by adding a point
$\si_*$ that projects to $q_0$.   The domain $\Si$ of $\si_*$ has $4$
components with $\Si_0,\Si_2,\Si_3$ all meeting  $\Si_1$ and a marked
point $z_0\in \Si_3$.
The map $h_0:\Si_0\to \De_J$ parametrizes $\De_J$ as a section, $h_1$ takes
$\Si_1$ onto the point $F_0\cap \De_J$, and  $h_2,h_3$ have image $F_0$
with $h_3(z_0) = x_0$. 
 The argument of Example~\ref{ex:lens} gives the following result.

\begin{lemma} The
space $\Vv_{2,J}$  of gluing parameters over $\Zz_{2,j}\cup\{\si_*\} =
S^2$ is the bundle $ \Oo(-2)\oplus \C$.  
\end{lemma}

\MS

\NI
{\bf \S 3.2.2 The stratum $\Zz_{1,J}$}\SS

This gives half of $\Ll_{\Zz}$.  The
other half comes from  the link of the orbifold $\Zz_{1,J}$ in $\Vv_{1,J}$.
Thus the next step is to look at $\Zz_{1,J}$.

Let $p_0, p_1$ be two distinct points on $F_0 \equiv S^2$, with $p_0=x_0$ 
and $p_1 = \Delta_J\cap F_0$.  Then $\Zz_{1,J}$ is the orbifold 
$$
\Zz_{1,J} = Y = \oMm_{0,2}(S^2,
p_0, p_1, 2)
$$
 of all stable maps to $S^2$ with two marked points $z_0, z_1$
that are in the class $2[S^2]$ and are such that $h(z_0) = p_0, h(z_1) = p_1$.
We will also need   to
consider the space 
$\oMm_{0,0}(S^2,2)$  of genus $0$ stable
maps of degree $2$ into $S^2$ that have no marked points and the space
$$ \TY= \oMm_{0,3}(S^2, p_0, p_1, p_2, 2)
$$
 of all degree $2$
stable maps to $S^2$ with three marked points $z_0, z_1, z_2$ such that
$h(z_0) = p_0, h(z_1) = p_1, h(z_2) = p_2 $.

\begin{lemma}\label{ss}  (i) $\oMm_{0,0}(S^2,2)$  is
a smooth manifold  diffeomorphic to $\C P^2$. \SS

\NI
(ii)  $\TY= \oMm_{0,3}(S^2, p_0, p_1, p_2, 2)$ is
a smooth manifold  diffeomorphic to $\C P^2$. \SS

\NI
(iii)  The forgetful map $f:\TY\to Y$ may be identified with the $2$-fold cover
that quotients $\C P^2$ by the involution $\tau:[x:y:z]\mapsto [x:y:-z].$ In
particular, $Y$ is smooth except at the point $\si_{01} = f([0:0:1])$
that has the local chart $\C^2/(x,y) = (-x,-y)$.  This point $\si_{01}$ is 
the stable map $[S^2, h, z_0, z_1]$ 
where  the critical values  of $h$ are at $p_0$ and $p_1$. 
\end{lemma}
\proof{} (i)  The space $\oMm_{0,0}(S^2,2)$ has
two strata. The first, $\Ss_1$, consists of self-maps of $S^2$ of degree $2$,
and the second, $\Ss_2$, consists of maps whose domain has two
components, each taken into $S^2$ by a map of degree $1$.  The equivalence
relation on each stratum is given by precomposition with a holomorphic
self-map of the domain.  lt is not hard to check that each equivalence class of
maps in $\Ss_1$ is uniquely determined by its two critical values (or
branch points).  Since these can be any pair of distinct points,  $\Ss_1$
 is diffeomorphic to the set of unordered pairs of distinct points in
$S^2$.  On the other hand there is one element $\si_w$ of $\Ss_2$  for each
point $w \in S^2$, the correspondence being given by taking $w$ to be the
image under $h$ of the point of intersection of the two components. 

If $\si_{\{x,y\}}$ denotes the element of $\Ss_1$ with critical values
$\{x,y\}$, we claim that $\si_{\{x,y\}} \to \si_w$ when $x,y$ both
converge to $w$. To see this, let
$h_{\{x,y\}}:S^2\to S^2$ be a representative  of   $\si_{\{x,y\}}$
and let $\al_{\{x,y\}}$ be the shortest geodesic from $x$ to $y$. (We assume
that $x,y$ are close to $w$.)  Then  $h_{\{x,y\}}^{-1}(\al_{\{x,y\}})$ is a circle
$\ga_{\{x,y\}}$  through the critical points of $h_{\{x,y\}}$.  This is obvious if
$h_{\{x,y\}}$ is chosen to have critical points at $0,\infty$ and if $x=
0,y=\infty$ since $h_{\{x,y\}}$  is then a map of  the form $z\mapsto
az^2$.  It follows in the general case because Mobius transformations take
circles to circles.   Hence $h_{\{x,y\}}$ takes each component of $S^2 -
\ga_{\{x,y\}}$ onto $S^2 - \al_{\{x,y\}}$.  If we now let $x,y$ converge to
$w$, we see that $\si_{\{x,y\}}$ converges to  $\si_w$.

The above argument shows that $\oMm_{0,0}(S^2,2)$ is the quotient of
$S^2\times S^2$ by the involution $(x,y)\mapsto (y,x)$.  This is well known
to be $\C P^2$.  In fact, it is easy to check that the map
$$
H: ([x_0:x_1],[y_0:y_1])\;\mapsto \; [x_0y_0: x_1y_1: x_0y_1 + x_1y_0 - 
x_0y_0- x_1y_1]
$$
induces a diffeomorphism from the quotient to $\C P^2$.  Under this
identification the stratum $\Ss_2 = H(diag)$ is the quadric
$(u+v+ w)^2 = 4uv$ (where we use coordinates $[u:v:w]$ on $\C P^2$). 
Further,  if we put
$$
p_0 = [0:1], \quad p_1 = [1:0],\quad p_2 = [1,1],
$$
the set of points in $\oMm_{0,0}(S^2,2) = \C P^2$ consisting of maps that
branch over $p_i$ is a line $\ell_i$, the image by $H$ of $(S^2\times p_i)\cup
p_i\times S^2$. Thus
$$
\ell_0 = \{u=0\},\quad \ell_1 = \{v=0\},\quad \ell_2 = \{w=0\}.
$$
Note finally that all stable maps in $\oMm_{0,0}(S^2,2)$ are invariant by an
involution: for example the map $z\mapsto z^2$ is invariant under the
reparametrization $z\mapsto -z$.  Since all elements have the same
reparametrization group, $\oMm_{0,0}(S^2,2)$ is smooth.  However, this will
no longer be the case when we add two marked points.\QED

\NI (ii)
Now consider
the forgetful map
$$
\phi_{30}: \oMm_{0,3}(S^2,p_0,p_1,p_2, 2) \to \oMm_{0,0}(S^2,2).
$$
For a general point of $\oMm_{0,0}(S^2,2)$, that is a
point where neither branching point is at $p_0$, $p_1$ or $p_2$, $\phi_{30}$
is $4$ to $1$.  To see this, note that for $i = 0,1,2$, $z_i$ can be either
of the points that get mapped to $p_i$ which seems to
give an $8$-fold cover.  However, because $h$ has degree $2$, $h$ is
invariant under an involution $\ga_h$ of $S^2$ that interchanges the two
inverse images of a generic point.   Hence the cover is $4$ to $1$, and the
covering group is  $\Z/2\Z \oplus \Z/2\Z$.  When just one branching point is
at some $p_i$, $\phi_{30}$ is $2$ to $1$, and when both branching points
are at some $p_i$, it is $1$ to $1$.   This determines $\phi_{30}$.  In fact, with
the above identification for $\oMm_{0,0}(S^2,2) = \C P^2$, $\phi_{30}$ is the
map $$
\phi_{30}: \C P^2\to\C P^2: \quad [x:y:z]\mapsto [x^2:y^2:z^2]. 
$$
Note  that the inverse image of $\Ss_2 = \{4uv=(u+v+w)^2\}$ consists of the
$4$ lines
$$
x\pm y \pm iz = 0.
$$
These components correspond to the $4$ different ways of arranging $3$
points on the two components of the stable maps in $\Ss_2$.  Note further
that none of the points in $\oMm_{0,3}(S^2,p_0,p_1,p_2, 2) $
are invariant by any reparametrization of their domains.  Hence all
points of this moduli space are smooth.\QED

\NI (iii)  Similar reasoning shows that the forgetful map
$$
\phi_{20}: Y=\oMm_{0,2}(S^2,p_0,p_1, 2) \to \oMm_{0,0}(S^2,2)
$$
 is a $2$-fold cover branched over $\ell_0\cup \ell_1$.  Hence we may
identify $Y$ as
$$
Y = \{[u:v:w:t]\in \C P^3: t^2 = uv\}
$$
where the cover $\phi_{20}:Y\to \C P^2$ forgets $t$.  There is one point in
$Y$ that is invariant under a reparametrization of its domain, namely the
point $\si_{01}$ corresponding to the map $h:S^2\to S^2$ that branches at
$p_0$ and $p_1$.  In the above coordinates on $Y$, 
$$
\si_{01} = \phi_2^{-1}(\ell_0\cap \ell_1) = [0:0:1:0].
$$
It is also easy to check that
$$
\phi_{32}: \oMm_{0,3}(S^2,p_0,p_1,p_2, 2) = \C P^2 \to Y
$$
has the formula
$$
\phi_{32}([x:y:z])= [x^2:y^2:z^2: xy].
$$
Since $\phi_{32}\circ \tau( [x:y:z]) = \phi_{32}( [x:y:-z]) = \phi_{32}([x:y:z])$,
$\phi_{32}$ is equivalent to quotienting out by $\tau$ as claimed.\QED

\MS

\NI
{\bf \S 3.2.3 The bundle  $\Vv_{1,J}\to \Zz_{1,J}$.}
\SS

Now we consider the structure of the orbibundles of gluing parameters over
$\TY = \oMm_{0,3}(S^2,p_0,p_1,p_3, 2)$ and $\Zz_{1,J}= Y=
\oMm_{0,2}(S^2,p_0,p_1, 2)$.  We will call the first $\TL\to \TY$ and the
second $L_Y\to Y$.  In both cases the fiber at the stable map $[\Si, h, z_i]$
is the tangent space  $T_{z_1}\Si$.

\begin{lemma} (i) The orbibundle $\TL\to \TY$ is smooth and may be
identified with the canonical line bundle $L_{can}$ over $\TY = \C P^2$.
\SS

\NI (ii) The orbibundle $L_Y\to Y$ is smooth except at the point
$\si_{01}$.   It can be identified with the quotient of $L_{can}$ by the
obvious lift $\tilde \tau$ of $\tau$.  
\SS

\NI(iii) The set $S(L_Y)$ of  unit vectors in $L_Y$ 
is smooth and diffeomorphic to $S^5$.  The orbibundle $S(L_Y)\to Y$
can be identified with the quotient of $S^5$ by the circle action
$$
\theta\cdot(x,y,z) = (e^{i\th} x, e^{i\th} y, e^{2i\th} z).
$$
\end{lemma}
\proof{}  Since $\TY$ is smooth, the general theory implies that $\TL$ is
smooth.  Therefore, it is a line bundle over $\C P^2$ and to understand its
structure we just have to figure out its restriction to one line.  It is easiest
 to consider one of the lines $x\pm y\pm iz=0$ that lie over
$\Ss_2$.  Recall that $\si_w\in \Ss_2$ is the stable map $[\Si_w,h_w]$
with domain $\Si_w = S^2\cup_{w=w}S^2$ and where $h$ is the identity map
on each component.  Suppose we look at the line in $\TY$ whose generic
point has $z_1$  on one component of $\Si_w$ and $z_0, z_2$ on the other. 
Then the bundle $\TL$ has a natural trivialization over the set $\{w\in S^2:
w\ne z_0,z_1,z_2\}$.  It is not hard to check that this trivialization extends
over the points $z_0, z_2$ but that one negative twist is introduced when
$z_1$  is added.  The argument is very similar to the proof of
Lemma~\ref{glu} below, and is left to the reader.\QED
\MS

\NI (ii) It follows from the general theory that $L_Y\to Y$ is smooth over the
smooth points of $Y$.  Moreover, at $\si_{01} = [S^2, h]$ the automorphism
$\ga:S^2\to S^2$ such that $h\circ \ga = h, \ga(z_i) = z_i$ acts on $T_{z_1}S^2$
by the map $v\mapsto -v$.  (To see this, note that we can identify $S^2$ with
$\C\cup \{\infty\}$ in such a way that $z_0 =p_0 = 0, z_1 =p_1 = \infty$. 
Then $h(z) = z^2$, and $\ga(z) = -z$.)  Hence the local structure of $L$ at
$\si_{01}$  is given by quotienting the trivial bundle $D^4\times \C$ by
the map $(x,y)\times v\mapsto (-x,-y)\times -v$.  This is precisely the
structure of the quotient of $\TL$ by $\tilde\tau$ at the singular point.  
Moreover, we can identify $S(L_Y)$  with $ S^5/\tau$ globally since
 $L_Y\to Y$ pulls back to $\TL\to \TY$ under the map $\TY\to Y$.\QED
\MS

\NI The quotient $S^5/\tau$ is
smooth except possibly at the fixed points $(x,y,0)$ of $\tau$.  Since
$S(L_Y)$ is smooth at these points, $S(L_Y)$ is smooth  everywhere.  It may
be identified with $S^5$ by the map
$$
S(L_Y) \equiv S^5/\tau \to S^5:\quad (x,y,z)\mapsto (x\,\sqrt{1+|z|^2},
y\,\sqrt{1+|z|^2}, z^2).
$$
The last statement may be proved by noting that the formula
$$
(x,y,z)\mapsto [x^2:y^2:z:xy]\in \C P^3
$$
defines a diffeomorphism from the orbit space of the given circle action to
$\C P^2/\tau = Y$. \QED
\MS

\NI
{\bf \S 3.2.4 Attaching the strata.}
\SS

The next step is to understand how the two strata $\Vv_{1,J}$ and
$\Vv_{2,J}$ fit together. The two zero sections $\Zz_{1,J}$ and $\Zz_{2,J}$
intersect at the point $\si_*$.  Recall that
the domain $\Si$ of $\si_*$ has $4$ components with
$\Si_0,\Si_2,\Si_3$ all meeting  $\Si_1$ and a marked point $z_0\in \Si_3$  . 
The map $h_0:\Si_0\to \De_J$ parametrizes $\De_J$ as a section, $h_1$ takes
$\Si_1$ onto the point $x_1=F_0\cap \De_J$, and  $h_2,h_3$ have image $F_0$
with $h_3(z_0) = x_0$.  The stratum of $\Zz_{1,J}$ containing $\si_*$ consists
just of this one point.   Hence the local coordinates of $\si_*$ in $\Zz_{1,J}$ are
given by  two gluing parameters $(a_0,a_1)$.  If we write $z_{ij}$ for the point
$\Si_i\cap \Si_j$, these are
 $$
 (a_0,a_1)\;\;
\mbox{where} \;\;a_0\in T_{z_{12}}\Si_1\otimes T_{z_{12}}\Si_2,\;\; 
a_1\in T_{z_{13}}\Si_1\otimes T_{z_{13}}\Si_3,\;\;.
$$
 Similarly, the local coordinates for a neighborhood of $\si_*$ in 
$\Vv_{1,J}$ are  $$
(b,a_0,a_1)
$$
where $(a_0,a_1)$ are as before and $b
\in T_{z_{01}}\Si_0\otimes T_{z_{01}}\Si_1$ is a gluing
parameter at the point $z_{01}$ where the component $\Si_0$
mapping to $\De_J$ is attached. On the other hand the natural coordinates for
a neighborhood of $\si_*$ in $\Vv_{2,J}$ are triples $(w,b,a)$ where  $b$ is a
gluing parameter at the point $z_{03}$ where the component $\Si_0$ that
maps to $\De_J$ is attached to the fixed fiber $\Si_3$, $w$ is the point where
the moving fiber $\Si_2$ (the one not containing $z_0$) is attached to $\Si_0$
and $a$ is a gluing parameter at $w$.

\begin{lemma}\label{glu} The attaching map $\al$ at $\si_*$ has the form
$(b,a_0,a_1)\mapsto (w_b, ba_0,  ba_1)$, where $b\ne 0$ and $\|b\|$ is
small.  Here the map $b\mapsto w_b$ identifies a small neighborhood of $0$
in $T_{z_{01}}\Si$ with a neighborhood
of $x_1$ in $\Delta_J$ in the obvious way.
\end{lemma}
\proof{}  
The attaching of $\Zz_{1,J}$ to $\Zz_{2,J}$ comes from  gluing  at the
point $z_{01}$ via the parameter $b$.   Thus we are gluing the ``ghost
component" $\Si_1$ to the component $\Si_0$ that maps to $\De_J$  in the
space of stable $A$-curves that are holomorphic for a {\em fixed} $J$.  (It is
only when one glues at $a_0$ or $a_1$ that one changes the homology class
of the curve $\De_J$ and hence has to change $J$.)   In particular, we can
forget the components $\Si_2, \Si_3$ of the domain $\Si$ of $\si_*$,
retaining only the points $z_{12}, z_{13}$ on $\Si_1$  where they are
attached.  Therefore we can consider the domain of the attaching map $\al$
to be the $2$-dimensional space  
$$  
\{[\Si_0\cup_{z_{01}} \Si_1,  z_{12}, z_{13}, h_\De;
b] :\quad b\in \C = T_{z_{01}}\Si_0\otimes T_{z_{01}}\Si_1\}, 
$$ 
and its range to be the space of all elements $[\Si_0, q_0, w,
h_\De, J]$ where $\Si_0 = S^2$ and $w$ moves in a small disc about $x_1$.
Here, the map $h_\De: \Si_0\to \De_J$  is fixed and 
parametrizes $\De_J$  as a section.  We can encode this by picking two points
$q_1, q_2$ in $\Si_0$ that are different from $q_0 = \pi_J(x_0)$ and
then considering $h_\De$ to be the map that takes these two marked
points to two other fixed points on $\De_J$.  Thus the attaching map
$\al$  is equivalent to the following map $\al'$ that attaches diffferent
strata in the moduli space $\oMm_{0,4}(S^2)$  of $4$ marked points  on
$S^2$: $$
\al': \{(\Si_0\cup\Si_1, q_1,q_2,z_{12}, z_{13}; b): b\in \C\} \to \{(\Si_0,
q_1,q_2,z_{12}, z_{13}) \in \oMm_{0,4}(S^2)\}.
$$
Here, each $\Si_i$ is a copy of $S^2$ as before.  On the left $q_1, q_2$ are
two marked points on $\Si_0$ and $z_{12}, z_{13}$ are two marked points on
$\Si_1$.  On the right, we should consider the three points $q_1, q_2,
z_{13}$ to be fixed, while $z_{12} = w$ moves, since this corresponds to
our previous trivialization of the neighborhood of $\si_*$ in $\Vv_{2,J}$.  
Thus $\al'$ may be considered as a map taking $b$ to $w_b = z_{12}\in
\De_J$.  

It remains to check that as $b$ moves once (positively)  around $0$,
$w_b$ moves once positively around $z_{13}$.  This follows by examining 
the identification of the glued domain 
$$
\Si_b = \left(\Si_0 - D(z_{01})\right)\cup_{gl_b} \left(\Si_1 -
D(z_{01})\right) $$
 with $\Si_0 = S^2$. Observe that the two points $q_1,q_2 $ in 
 $\Si_0 - D(z_{01})$  and the single 
point $z_{12}$ in $\Si_1 -
D(z_{01})$ must be taken to the corresponding three fixed points on
$S^2=\Si_1$.  Hence the identification on $\Si_0 - D(z_{01})$ is fixed, while
that on $\Si_1 - D(z_{01})$ can rotate about $z_{13}$ as $b$ moves.  Hence,
when $b$ moves round a complete circle, so does $w_b$.  It remains to
check the direction of the rotation.  Now, as we saw in 
Proposition~\ref{nbhdpt}, as $b$ moves once  round this circle positively  as
seen from $z_{01}$, the point $p_b$ on $\p D(z_{01})\subset \Si_1$ that is
matched with a fixed point $p$ on $\Si_0 -  D(z_{01})$ moves once positively
round $\p D(z_{01})$.  In order to line up $p_b$ with $p$, 
$\Si_1$ must be rotated in the {\em opposite} direction, i.e. positively as seen
from the fixed point $z_{13}$.  Hence $w_b$ rotates positively round
$z_{13}$.

To complete the proof of the lemma, we must understand how the gluing
parameters $a_0, a_1$ fit into this picture.  Since nothing is happening in
the vertical (ie fiberwise) direction, we may consider the $a_i$ to be
elements of the following tangent spaces:
$$
a_0\in T_{z_{12}} \Si_1,\quad a_1\in T_{z_{13}} \Si_1.
$$
As $b$ rotates positively, the image of $a_0$ in the glued curve rotates once
positively in the tangent space of $z_{12}$, and $a_1\in T_{w_b}\Si_1$ also
rotates once with respect to the standard trivialization of the tangent spaces
$T_{w_b}\Si_1 \subset T(\Si_1)|_{D(z_{12})}.$
Hence result.\QED

\NI
{\bf Proof of Proposition~\ref{prop:link}}.

We have identified the orbibundle $\Vv_{1,J}\to \Zz_{1,J}$ with $L\to Y$
 and the bundle $\Vv_{2,J}\to \Zz_{2,J}$ with $\Oo(-2)\oplus
\C\to S^2$. The
previous lemma shows that these are attached by first twisting $\Vv_{1,J}$
to $\Oo(-3)\oplus \Oo(-1)$ and then plumbing it to $L$.  Hence
$$
\Ll_\Zz = S(\Oo(-3)\oplus \Oo(-1))\pl S(L_Y)
$$
as claimed.     The identification of the latter space with
$S(\Oo(-1)\oplus\C)$ follows from
Proposition~\ref{str}.         
\QED

\subsection{The projection $\Vv_J\to \Jj$}

In order to complete the calculation of the link $\Ll_{2,0}$ of $\Jj_2$ in $\Jj$
it remains to understand  the projection $\Vv_J\to \Jj$.  This is $1$-to-$1$
except over the points of $\Jj_1$.  In $\Vv_{2,J}$ it is clearly the points 
with zero gluing parameter at the moving fiber that get collapsed.  Thus the
subbundle $R_-$ of the circle bundle 
$S(L_P) \to \Pp(\Oo(-2)\oplus \C)$ that lies over 
the (rigid) section $S_-= \Pp(\{0\}\oplus \C)$ must be collapsed to a single
circle. The subbundle $R_+$ lying over the other section  
$S_+=\Pp(\Oo(-2)\oplus \{0\})$  maps to a family of distinct elements in
$\Jj_1$.  

The story on $\Vv_{1,J}$ is, of course, more complicated.  Here the points
that concern us are the maps in $\Ss_2$
where the branch points coincide.  Thus, if we identify
$\oMm_{0,0}(S^2, 2)$ with $\C P^2$ as in Lemma~\ref{ss}, these are the
points of the quadric $Q = \{(u+v+w)^2 = 4uv\}$.   
Note that the attaching point $\si_*\in \oMm_{0,2}(S^2, p_0,p_1, 2)$ sits over
$$
[1:0:-1]=\ell_1 \cap
Q\;\in\;\C P^2 = \oMm_{0,0}(S^2, 2).
$$

The lift of $Q$ to $ \oMm_{0,2}(S^2, 2)$ has two components $Q_\pm$, 
given by the intersection  $Y\cap H_\pm$ where $H_\pm$ is the
hyperplane $2t = \pm (u+v+w)$. Since we can assign these at will, we will say
that $Q_-$ corresponds to elements with the two marked points $z_0, z_1$
on the same component of $\Si_w = \Si_0\cup_{w=w} \Si_1$ and that $Q_+$
corresponds to elements with $z_0, z_1$ on different components.  Then,
when one glues at $z_1$ the resulting $A$-curve is the union of an
$(A-F)$-curve  with an $F$ curve.  It is not hard to check that the points on
$Q_-$ give rise to a $(A-F)$-curve through $x_0$, which is independent of
$w$, while those on $Q_+$ give rise to a varying $(A-F)$-curve that meets
the $J_w$-holomorphic fiber through $x_0$ at the point $w$.  Note that the
intersection
 $Q_+\cap Q_-$ consists of two points, $p_* = [1:0:-1:0]$ (corresponding to
$\si_*$) and $q_*=[0:1:-1:0]$.
 Moreover, in the coordinates
$(a_0, a_1)$ of a neighborhood of $\si_*$ used in Lemma~\ref{glu} above, 
$$
 \{(a_0, a_1): a_0=0\}\subset Q_-, \quad  \{(a_0,a_1): a_1=0\}\subset Q_+.
$$
This confirms that when $\Vv_{2,J}\otimes \Oo(-1) = \Oo(-3)\oplus
\Oo(-1)$ is plumbed to $\Vv_{1,J}$, $Q_-$ is plumbed to the subbundle
$\{0\} \oplus \Oo(-1)$ corresponding to $R_-$ and $Q_+$ is plumbed to
the subbundle  $\Oo(-3)\oplus\{0\}$ corresponding to $R_+$.

Let $S(Q_\pm)\to Q_\pm$ 
 denote the restriction of $S(L_Y)\to Y$ to $Q_\pm.$ Then the plumbing
$\Oo(-3)\oplus\Oo(-1) \pl S(L_Y)$ contains the plumbings $R_-\pl
S(Q_-),$ and $ R_+\pl S(Q_+).$ 

\begin{lemma}  (i) $R_-\pl
S(Q_-) = S(\C) =S^2\times S^1$ and $ R_+\pl S(Q_+) = S(\Oo(-2).$ 
\SS

\NI(ii)  The subsets $R_-\pl S(Q_-)$ and $R_+\pl S(Q_+)$ of 
$\Oo(-3)\oplus\Oo(-1) \pl S(L_Y)$
intersect in a circle.
\end{lemma}
\proof{}  Since $Q_-$ and $Q_+$ do not meet the singular point of $Y$, both
bundles $S(Q_\pm)\to Q_\pm$ have Euler number $-1$.  Hence $$
R_-\pl  S(Q_-) = S(\Oo(-1))\pl S(\Oo(-1)) = S(\C) = S^2\times S^1,
$$
and
$$
R_+\pl S(Q_+) = S(\Oo(-3))\pl S(\Oo(-1)) = S(\Oo(-2)).
$$
This proves (i).  To prove (ii) note that the inverse image (in $S(Q_\pm)$) of
the intersection point $p_* = [1:0:-1]$ of $Q_-$ with $Q_+$ disappears under
the plumbing.  But the other one remains.\QED

\NI
{\bf Proof of Theorem~\ref{LINK}}.

It follows from part (i) of the preceding lemma that it is possible to collapse
 the subset $R_-\pl S(Q_-)$ of $\Ll_Zz $
to a single circle.  Moreover, it is not hard to see that 
 under the identification of   $\Ll_Zz=S(\Oo(-3)\oplus \Oo(-1))\pl S(L_Y)$
 with $S(\Oo(-1)\oplus\C)$, 
this collapsing corresponds to collapsing the circle bundle over the exceptional
divisor.   Since the intersection of $R_-\pl S(Q_-)$ with $R_+\pl S(Q_+)$ is a
single circle,  this collapsing does not affect $R_+\pl
S(Q_+)$.  Note that  $R_+\pl
S(Q_+)$ is the inverse image of some $2$-dimensional submanifold of $\C
P^2$.   Because $ R_+\pl S(Q_+) = S(\Oo(-2)$ this submanifold must be a
quadric.  \QED

\section{Analytic arguments}

In \S 4.1 we prove the (easy) 
Lemmas~\ref{mfld}.  \S 4.2 contains a detailed analysis
of gluing.  The exposition here is fairly
 self-contained, though some results are quoted from [MS] and [FO].

\subsection{Regularity in dimension $4$}

The theory of $J$-holomorphic spheres in dimension $4$ is much simplified
by the  fact that  any
$J$-holomorphic map $h:S^2\to X$ that represents a class $A$ with $A\cdot
A\ge -1$ is regular, i.e. the linearized delbar operator 
$$
Dh: W^{1,p}(h^*(TX))\to L^p\left(\La_J^{0,1}(S^2)\otimes h^*(TX)\right)
$$
 is surjective.  This remains true even if $h$ is a multiple covering.  (For a
proof see Hofer--Lizan--Sikorav~[HLS].  The notation is explained
in \S3.1 below.)  Therefore, regularity is automatic:
one does not have to perturb the equation in order to achieve it. The
analogous  statement when $A\cdot A < -1$ is that $\Coker Dh$ always has
 rank equal to $2 +2A\cdot A$.  As is shown below, this almost immediately
implies that  the $\Jj_k$ are 
submanifolds of $\Jj$.
\MS

 \NI
{\bf Proof of Lemma~\ref{mfld}}

 We begin by proving that $\Jj_k$ is a Fr\'echet
manifold.  This is obvious when $k = 0$, since $\Jj_0$ is an open subset of
$\Jj$. For $k > 0$, let $\Cc_k$ denote the space of all symplectically
embedded spheres in the class $A-kF$,
and let $\Cc_k(\Jj)$ be the bundle over $\Cc_k$ whose fiber at $C$ is the
space of all smooth almost complex structures on $C$ that are compatible with
$\om|_C$.  Then $\Cc_k(\Jj)$ fibers over $\Cc_k$ and it is easy to check
that both spaces are Fr\'echet manifolds.  (Note that $\Cc_k$ is an open
submanifold in the space of all embedded spheres in the class $A-kF$.
Because these spheres are not parametrized the tangent space to $\Cc_k$ at
$C$ is the space of all sections of the normal bundle to $C$.) Further $\Jj_k$
fibers over $\Cc_k(\Jj)$ with fiber at $(C,J|_C)$ equal to all $\om$-compatible
almost complex structures that restrict to $J$ on $TC$.  This proves the claim.

To see that $\pi_k$ is bijective when $k > 0$ note that each $J\in \Jj_k$
admits a holomorphic curve in class $A - kF$ by definition, and that this
curve is unique by positivity of intersections.  A similar argument works
when $k = 0$ since the curves in $\Mm(A,\Jj)$ are constrained to go through
$x_0$. Hence $\Mm_k$ inherits a Fr\'echet manifold structure from $\Jj_k$.

To show that $\Jj_k$ is a submanifold of $\Jj$ when $k > 0$ we must 
use the theory of $J$-holomorphic curves, as explained in Chapter~3 of [MS]
for example.
Let $\Mm_k^s, \Jj_k^s, \Jj^s$ denote the similar spaces in the
$C^s$-category for some large $s$.  These are all Banach manifolds. 
It is easy to check that the tangent space $T_J\Jj^s$ is the space
$End(TX,\om,J)$ of all $C^s$-sections $Y$ of the endomorphism bundle of
$TX$ such that $$
JY + YJ = 0,\quad \om(Yx,y) = \om(x, Yy).
$$
These conditions imply that $\om(Yx,x) = \om(Yx,Jx) = 0$ for all $x$.  It
follows easily that $Y$ is determined by its value on a single nonzero
vector $x$ that it has to take to the $\om$-orthogonal complement to the
$J$-complex line through $x$.
Observe further that there is an exponential map
$$
exp: T_J\Jj^s\to \Jj^s
$$
that  preserves
smoothness and is a local diffeomorphism near the zero section.

Next, note that the tangent space $T_{[h,J]}\Mm_k^s$ is the quotient of
the space of all pairs $(\xi, Y)$ such that
$$
Dh (\xi) + \frac 12 Y\circ dh\circ j = 0\qquad\qquad (*)
$$
by the $6$-dimensional tangent space to the reparametrization group
$\PSL(2,\C)$.  Here $j$ is the standard almost complex structure on $S^2$
and $Dh$ is the linearization of the delbar operator that maps 
the Sobolev space of $W^{1,p}$-smooth sections of
$h^*(TX)$ to anti-$J$-holomorphic $1$-forms, viz:
\begin{eqnarray}\label{eq:Dh}
Dh: W^{1,p}(S^2, h^*(TX)) \to L^p(\La_J^{0,1}(S^2, h^*(TX)),
\end{eqnarray}
where the norms are defined using the standard metric on $S^2$ and a
metric on $TX$.

In~[HLS], Hofer--Lizan--Sikorav show
how to interpret elements of $\Ker\, Dh$  and of $\Ker\, Dh^*$ (where $Dh^*$ is
the formal adjoint) as $J$-holomorphic curves in their own right.  Using the fact
that the domain is a sphere and that $X$ has dimension $4$, they then use
positivity of intersections to show that $\Ker\, Dh$ is trivial when $k > 0$, i.e. it
consists only of vectors that generate the action of $\PSL(2,\C)$.  Hence $\Ker\,
Dh^*$ is a bundle over $\Mm_k^s \cong \Jj_k^s$ of rank $4k - 2 = -{\rm
index\,}Dh$, and it is not hard to see that it is isomorphic to the normal
bundle of $\Jj_k^s$ in $\Jj^s$.  In other words 
$$
T_J\Jj^s = T_J \Jj_k^s \oplus \Ker\, Dh^*.
$$
To see this, observe that the map 
\begin{eqnarray}\label{eq:io}
\io: Y\mapsto \frac 12 Y\circ dh\circ j
\end{eqnarray}
maps $T_J\Jj^s$
onto the space of $C^s$-sections of $\La_J^{0,1}(S^2, h^*(TX))$, and that the
kernel of this projection consists of elements $Y$ that vanish on
the tangent bundle to the image of $h$
 and so lie in $T_J \Jj_k^s$ whenever $[h,J] \in \Mm_k^s$.
It follows from equation~(\ref{eq:Dh}) above that the image of $T_J \Jj_k^s$
under this projection is precisely equal to the image of $Dh$, and so has
complement isomorphic to $\Ker\, Dh^*$.  (For more details on all this, see the
Appendix to~[A].)   

It now remains to show that $\Jj_k$ is a submanifold of $\Jj$
whose normal bundle has fibers $\Ker\, Dh^*$.  This means in particular
that the codimension of $\Jj_k$ is $- \ind\,Dh = 4k-2.$  We  therefore
have to check that each point in $\Jj_k$ has a neighborhood $U$ in $\Jj$ that
is diffeomorphic to the product $(U\cap \Jj_k)\times \R^{4k-2}$.  It is here
that we use the exponential map {\it exp}.  Clearly, one can use {\it exp} to define 
such local charts for $\Jj_k^s$ in $\Jj^s$.  The point here is that the
derivative of the putative chart will be the identity along  $(U\cap
\Jj_k^s)\times \{0\}$ and so by the implicit function theorem for Banach
manifolds will be a diffeomorphism on a neighborhood.  Then, because
$\Ker\, Dh^*$ consists of $C^\infty$ sections when $J$ is $C^\infty$ and
because {\it exp} respects smoothness, this local diffeomorphism will take 
$(U\cap \Jj_k)\times \R^{4k-2}$ onto a neighborhood of $J$ in $\Jj$.
\QED

\subsection{Gluing}

The next task is  to complete the proof of Propositions~\ref{nbhdpt}
and~\ref{glue}.
The standard gluing methods are local and work in the neighborhood of one
stable map, and so our main problem is to globalize the construction. 
The first step in doing this is to show that one can still glue even when the
elements of the obstruction bundle are nonzero at the gluing point. 
We will use the gluing method of McDuff--Salamon~[MS]
and Fukaya--Ono~[FO].  Much of the needed analysis appears in~[MS] 
 but the conceptual framework of that work has to be
enlarged to include the idea of stable maps as in Hofer--Salaomon~[HS].  No
doubt the other gluing methods can be adapted to give the same results.  

Our aim is to construct a gluing map
$$
\Gg:\Nn_\Vv(\Zz)\to \oMm(A-kF,\Jj)
$$
where $\Zz = \oMm(A-kF, \Jj_m)$  is the space of stable maps in
class $A-kF$ with one component in class $A-mF$, and 
$\Nn_\Vv(\Zz)$ is a neighborhood
of  $\Zz$ in the space $\Vv$ of gluing parameters.  Choose once and for all a 
$(4m-2)$-dimensional
subbundle $K$ of $T\Jj|_{\Jj_m}$ that is transverse to $\Jj_m$. 
 As explained in \S3.1 above the exponential map $exp$ maps a
neighborhood
 of the zero section in $K$ diffeomorphically onto a
neighborhood of $\Jj_m$ in $\Jj$.   For each $J\in \Jj_m$ let 
$$
\Kk_J\subset\Jj
$$
 be the slice  through $J$ (i.e. the image under  $exp$ of a small
neighborhood $\Nn_J(K)$ of $0$ in the fiber of $K$ at $J$).  We shall prove the
following sharper version of  Proposition~\ref{glue}.

\begin{prop}\label{glue2}  Fix $J\in \Jj_m$ and let $\Nn_\Vv(\Zz_J)$ be the
fiber of the map $\Nn_\Vv(\Zz)\to \Jj_m$ at $J$.  Then, if the
neighborhood $\Nn_\Vv(\Zz)$ is sufficiently small, there is a 
homeomorphism  
$$
\Gg_J:\Nn_\Vv(\Zz_J)\;\longrightarrow\; \oMm(A-kF,\Kk_J)
$$
 onto a neighborhood $\oMm(A-kF, J)$ in $\oMm(A-kF,\Kk_J)$. Moreover,
the union of all the sets $\im \Gg_J, J\in \Jj_m$, is a neighborhood of
$\oMm(A-kF,\Jj_m)$ in $\oMm(A-kF,\Jj)$. \end{prop}

Let $\pi_J:\Nn_\Vv(\Zz_J)\to\Zz_J$ denote the projection.  
We will first construct the map $\Gg_J$ in the fiber
at one point $\si = [\Si_\si, h_\si, J]$ of $\Zz_J$ and then how to fit these
maps together to get a global map over  $\Nn_\Vv(\Zz_J)$ with the stated
properties.  For the next paragraphs (until \S4.2.4) we will fix a particular
representative $h_\si:\Si_\si\to X$ of $\si$, and we will define $\TGg$ as a map
into the space of parametrized stable maps.  In order to understand a full
neighborhood of $\si$ we will have to glue not only at points where the
branches meet the stem $\Si_0$ but also at points internal to the branches. 
Therefore, for the moment we will forget the stem-branch structure of our
stable maps and consider the general  problem of gluing, at the points
$z_i\in\Si_{i0}\cap\Si_{i1}$ with parameter 
$$
a\; = \;\oplus_i  a_i\;\in \; \bigoplus_i T_{z_i}\Si_{i0}\otimes
T_{z_i}\Si_{i1}. $$
 \MS

\NI
{\bf  \S 4.2.1:\,  Construction of the pregluing $h_a$}\SS

We showed in Proposition~\ref{glue} above
 how to construct the glued domain $\Si_a$.  Since
this construction depends on a choice of metric on $\Si$, we must assume that
the domain $\Si$ of each stable map is equipped with a K\"ahler metric that
is flat near all double points and is invariant under the action
of the isotropy group $\Ga_{\si}$. 
Fukaya--Ono point out in [FO] \S 9 that it is possible to choose such a
metric continuously over the whole moduli space: one just has to start at
the strata containing elements $\si$ with the largest number of
components, extend the choice of metric near these strata by using the
gluing construction  (which is invariant
by $\Ga_\si$) and then continue inductively, strata by strata.   In what
follows we will assume this has been done.  We will also suppose that the
cutoff functions $\chi_r$ used to define $\Si_a$ have been chosen once
and for all.

The approximately holomorphic map $h_a:\Si_a\to X$ is defined from
$h_\si$ by using cutoff functions.  As before, we write  $r_i$ or simply $r$
instead of $\sqrt{\|a_i\|}$.  Hence if $R$ is as in [FO] or [MS], $r= 1/R$.  We
choose a small $\de> 0$ once and for all so that $r/\de $ is still small.\footnote
{
The logic is that one chooses $\de>0$ small enough for certain inequalities to
hold, and then chooses $r\le r(\de)$.  See Lemma~\ref{le:de} below.}
 Set $x_i = h_\si(z_i)$.   Then, for $\al = 0,1$
define  
\begin{eqnarray*}
h_a (z) &=& h_\si(z) \mbox{ for }\; z\in \Si_{i\al} - D_{z_i}^\al(2r/\de)\\ 
& =
&   x_i  \mbox{ for }\; z\in D_{z_i}^\al({r}/{\de } ) -
D_{z_i}^\al(r) \end{eqnarray*}
 and interpolate 
 on the annulus $D_{z_i}^\al(2{r}/{\de} ) - D_{z_i}^\al(r/\de)$ in $\Si_{i\al}$
 by setting
$$
h_a(z) = exp_{x_i}(\rho( {\de |z|}/ r) \xi_{i\al}(z)),
 $$
where  $\rho$ is a smooth cut-off function that equals $1$ on $[2,\infty)$
and $0$ on $[0,1]$, and
 the vectors $\xi_{i\al}(z)\in T_{x_i} X$ exponentiate to give $h_\si(z)$
on $\Si_{i\al}$:
$$
h_\si(z) = exp_{x_i}(\xi_{i\al}(z)),\mbox{ for }\; z\in 
D_{z_i}^\al( 2r/{\de}). $$
The whole expression is defined provided that $2r/\de$
is small enough for the exponential maps to be injective.

Later it will be useful to consider the corresponding map $h_{\si,r}$ with
domain $\Si$.   This map equals $h_a$ on  $\Si - \cup_{i,\al} D_{z_i}^\al
(r_i)$ and is set equal to $x_i$ on each disc $D_{z_i}^\al(r_i)$.  Note that
$h_{\si,r}:\Si\to X $ converges in the $W^{1,p}$-norm to $h_\si$ as
$r\to 0$.  

\MS

\NI{\bf \S 4.2.2\, Construction of the gluing $\TGg(h_\si,a)$.}\SS

Let 
$$
\Nn_0(W_a) = \Nn_0((W^{1,p}(\Si_a, h_a^*(TX)))
$$
 be a small
neighborhood of $0$ in $W^{1,p}(\Si_a, h_a^*(TX)))$.  Note that, if $\Si_a$
has several components $\Si_{a,j}$, the elements $\si$ of $W_a$ can
be considered as collections $\xi_j$ of  sections in $W^{1,p}(\Si_{a,j},
(h_{a,j})^*(TX))$  that agree pairwise at the points $z_i$. 
(This makes sense since the $\xi_j$ are continuous.) Further, we may
identify $\Nn_0(W_a)$
 via the exponential
map with a neighborhood of $h_a$ in the space of $W^{1,p}$-maps $\Si_a\to
X$. We will write $h_{a,\xi}$ for the map $\Si_a\to X$ given by:
$$
h_{a,\xi}(z) = exp_{h_a(z)}(\xi(z)), \quad z\in \Si_a.
$$

Recall that $\Nn_J(K)$ is a neighborhood of $0$ in the fiber of $K$. 
Given $Y\in \Nn_J(K)$ we
will write $J_Y$ for the almost complex structure $exp(Y)$ in the slice $
\Kk_J$.
Now consider
the locally trivial bundle $\Ee = \Ee_a \to \Nn_0(W_a)\times \Nn_J(K)$ whose
fiber at $(\xi, Y)$ is  
$$
\Ee_{(\xi,Y)} = L^p(\La^{0,1}(\Si_a)\otimes_{J_Y} h_{a,\xi}^*(TX)). 
$$
We wish to convert the pregluing $h_a$ to a map that is $J_Y$-holomorphic 
for some $Y$  by
using the implicit function theorem for the section $\Ff_a $ of $\Ee_a$
defined by  $$
\Ff_a(\xi, Y) = \op_{J_Y}(h_{a,\xi}).
$$
 Note that
$\Ff_a(\xi,Y) = 0$ exactly when the map $h_{a,\xi}$ 
is $J_Y$-holomorphic.

 The linearization $\Ll(\Ff_a)$ of $\Ff_a$
at  $(0,0)$ equals 
$$
\Ll(\Ff_a) = D(h_a)\oplus \io_a: W^{1,p}(h_a^*(TX))
\oplus 
K \;\to \;
L^p(\La^{0,1}(\Si_a)\otimes_{J} h_a^*(TX)),
$$
where $\io_a$ is defined by $\io_a(Y) = \frac 12 Y\circ
dh_a\circ j$ as in equation~(\ref{eq:io}) in \S4.1.  

\begin{lemma}\label{bound}  Suppose that there is a 
continuous family of right inverses
$Q_a$ to $\Ll(\Ff_a)$ that are uniformly bounded 
 for $\|a\|\le r_0$. Then,
there is $r_1 > 0$ such that for all $a$ satisfying $\|a\|\le r_1$ there is a
unique element $ (\xi_a,Y_a)\in \im Q_a$ such that  
$$
\Ff_a(\xi_a, Y_a)=0.
$$
Moreover, $(\xi_a,Y_a)$ depends continuously on the initial data.
\end{lemma} 
\proof{}  This follows from the implicit function
theorem as stated in 3.3.4 of [MS].  It also uses Lemma A.4.3 of [MS].  See also
[FO] \S11.\QED

We will construct the required family $Q_a$ in \S 4.2.3.  
By the above lemma, this allows us to define the gluing map.

\begin{defn}\label{defgl1}\rm  We set $\TGg(h_\si, a) = (\Si_a, h_{a,\xi_a},
J_{Y_a})$ where $(\xi_a,Y_a)$ is the unique element in the above lemma. 
Further  $\Gg(h_\si,a)] = [\Si_a, h_{a,\xi_a},
J_{Y_a}]$.
\end{defn}

The next proposition states the main local properties of the gluing map $\Gg$.

\begin{prop}\label{exun}  Each $\si\in \Zz_J$ has a neighborhood
$\Nn_\Vv(\si)$ in $\Vv_J$ such that the map
$$
\Nn_\Vv(\si)\to \oMm(A-kF,\Kk_J): \quad (\si',a') \mapsto \Gg(h_{\si'}, a')
$$
 takes $\Nn_\Vv(\si)$ bijectively onto an open subset in $
\oMm(A-kF,\Kk_J)$.  Moreover this map depends continuously
on $J \in \Jj_m$.
\end{prop} 
\proof{}  This is a restatement of Theorem
12.9 in [FO].   Note that the stable map
$\Gg(h_{\si'}, a')$  depends on the choice of representative $(\Si', h_{\si'})$ of
the equivalence class $\si' = [\Si', h_{\si'}]$.  However, it is always possible to
choose a smooth family of such representatives in  a small enough
neighborhood of $\si$ in $\Zz_J$.  (This point is discussed further in \S 4.2.4.)  
Moreover, if $\si$ is an orbifold point (i.e. if $\Ga_\si$ is nontrivial), then
$h_\si$ is $\Ga_\si$-invariant and one can define  $\TGg$ so that it is
equivariant with respect to the natural action of $\Ga_\si$ on the space of
gluing parameters $a$ and its action on a neighborhood of $\si$ in the space of
parametrized maps.  The composite $\Gg$  of $\TGg$ with the forgetful map is
therefore
$\Ga_\si$-invariant.
 (Cf. the discussion before Lemma~\ref{le:repr}.)  This
shows that $\Gg$ is well defined.

One proves that it is a local homeomorphism as in [FO] \S 13, 14, and we will
say no more about this except to observe  that our adding of
$K$ to the domain of $Dh_\si$ is equivalent to their replacement of the
range of $Dh_\si$ by the quotient $L^p/ \io_a(K)$.  \QED

\NI
{\bf \S 4.2.3\, Construction of the right inverses $Q_a$}\SS

This is done essentially as in A.4 of [MS] and \S12 of  [FO].  However, there
are one or two extra points to take care of, firstly because the stem of the
map $h_\si$   is not regular, so that the restriction of $D{h_\si}$ to $\Si_0$
is not surjective, and secondly because the elements of  the normal
bundle $K\to \Jj_m$ 
 do not necessarily vanish near the  points
$x_i$ in $X$ where gluing takes place.  

For simplicity, let us first consider the case when  $\Si$ has just two
components $\Si_0, \Si_1$  intersecting at the point $w$, and that $h_{\si}$
maps $\Si_0$ onto the $(A-kJ)$-curve $\De_J$ and   $\Si_1$ onto a fiber.  (For
the general case see Remark~\ref{rmk:gen}.)  
Then the linearization of $\op_J$ at $h_\si$ 
 has the form   
$$
D{h_\si}: W^{1,p}(\Si, h_\si^*(TX)) \to L^p(\La_J^{0,1}(\Si)\otimes
h_\si^*(TX)). 
$$
Here the domain consists of pairs $(\xi_0,\xi_1)$, where $\xi_j$ is a
$W^{1,p}$-smooth section of the bundle $h_{\si_j}^*(TX)\to\Si_j$, subject to
the condition $$
\xi_0(w) = \xi_1(w),
$$
 and the range consists of
pairs of $L^p$-smooth $(0,1)$-forms over $\Si_j$ with values in 
$h_{\si_j}^*(TX)$ and with no
condition at $w$.  For short we denote this map
by $$
Dh_\si: W_\si \to L_{\si_0}\oplus L_{\si_1}.
$$

Recall from the discussion before Proposition~\ref{glue2}, that we chose
$K$ so that  
$$
Dh_{\si_0}\oplus \io_0: W_{\si_0}\oplus K \to L_{\si_0}.
$$
is surjective and $
\io_0: K\to  L_{\si_0}
$
is injective.  (All maps $\io$ are defined as in equation~(\ref{eq:io}): it should
be clear from the context what the subscripts mean.)

\begin{lemma}\label{le:NJ}  There are constants  $c,r_0 > 0$ so that the
following conditions hold for all $r < r_0$: \SS

\NI
(i) $\io_a$ is injective for all $\|a\| \le r$;

\NI
(ii) the projection  $pr_K: L_{\si_0} \to K$ that has kernel $\im
D_{\si_0}$ and satisfies $pr_K\circ \io_0 = {\rm id}_K$ has norm $\le c$;

\NI
(iii)  for all $Y\in K$ and $j = 0,1$
$$
\left(\int_{D_{w}^j(r)} |\io_{j} (Y)|^p\right)^{1/p}  \le \frac 1{12 c}
\left( \int_{\Si_j}
 |\io_j (Y)|^p\right)^{1/p},
$$
where ${D_{w}^j(r)}$ is the disc in $\Si_j$ on which gluing takes place
and  integration is with respect to the area form defined by the chosen
K\"ahler metric on $\Si_\si$. 
\end{lemma}
\proof{}   There is $c$ so that (ii) holds
because $\im D_{\si_0}$ is closed and $\im\io_0$ is finite
dimensional.  Then there is $r_0 = r_0(c)$ satisfying (i) and (iii)  since  the
elements of $K$ are $C^\infty$-smooth (as are the elements of $\Jj$).  \QED

\begin{lemma}\label{dh}  The operator
$$
Dh_\si\oplus (\io_0, \io_1): W_\si\oplus K \to L_{\si_0}\oplus L_{\si_1},
$$
is surjective and has kernel ${\rm ker\,}Dh_\si$.
\end{lemma}
\proof{}   We know from the proof of Lemma~\ref{mfld} that
$$
Dh_{\si_0}\oplus \io_0: W^{1,p}(\Si_0, h_{\si_0}^*(TX))\oplus K \to
L^p(\La^{0,1}(\Si_0)\otimes_{J} h_{\si_0}^*(TX)) =L_{\si_0}, 
$$
is surjective.  Similarly, $Dh_{\si_1}$ is surjective.  Therefore, to prove
surjectivity we just need to check that the compatibility condition $\xi_0(w) =
\xi_1(w)$ for the elements of $W_\si$ causes no problem. 
However, the pullback bundle
$h_{\si_1}^*TX$ splits naturally into the sum of a line bundle with Chern class
$2d$ (where $d\ge 0$ is the multiplicity of $h_{\si,1}$) and a trivial line
bundle, the pullback of the normal bundle to the fiber $\im h_{\si_1}$. 
Hence there is a element $\xi_1$ of $\ker Dh_{\si_1}$ with any given
value $\xi_1(w)$ at $w$.  The result follows.  Note that an appropriate
version of  this argument applies for all $\si$, not just those with two
components, since there is just one condition to satisfy at each
double point $z$ of $\Si$ and the  maps
$\ker Dh_{\si_j}\to \C^2: \xi\mapsto \xi(z)$ are surjective for $j > 0$.
The second statement holds because $\io_0$ is injective.
\QED  

Note that  the
right inverse $Q_\si$  to   $Dh_{\si}\oplus (\io_0,
\io_1)$  is completely determined by choosing a complement to
the finite dimensional subspace $\ker Dh_\si$ in $W_\si$.  
Consider the composite
$$
pr_{\si_0}: L_{\si_0} \oplus L_{\si_1}\longrightarrow L_{\si_0}
{\longrightarrow} K 
$$
where the second projection is as in Lemma~\ref{le:NJ} (ii).
The fiber $(pr_{\si_0})^{-1}(Y)$ at $Y$ has the form $(\im Dh_{\si_0} +
\io_0(Y))\oplus L_{\si_1}$, and we write
$$
Q_\si^Y: \im Dh_{\si_0}\oplus L_{\si_1}  + (\io_0(Y), \io_1(Y)) \to W_\si
$$
for the restriction of  $Q_\si$  to this fiber.

We now use the method of [MS] A.4 to construct an approximate  right
inverse $Q_{a,app}$  to  $$
\Ll(\Ff_a) = Dh_a\oplus \io_a: W_a
\oplus 
K \;\to \; L_a,
$$
where
$$
W_a = W^{1,p}(h_a^*(TX)),\qquad L_a = L^p(\La^{0,1}(\Si_a)\otimes_{J}
h_a^*(TX)).
$$
It will be convenient  to use the
approximations  $
h_{\si,r}:\Si_\si\to X
$ to $h_\si$
that were defined at the end of \S 4.2.1 where $r^2 = \|a\|$.  We write
$h_{\si_j, r}$ for the restriction of $h_{\si, r}$ to the component $\Si_j$.  Since
$h_{\si,r}$ converges $W^{1,p}$ to $h_\si$ as $r \to 0$, $Dh_{\si,r}$ has  a
uniformly bounded inverse 
$$ 
Q_{\si,r}:\quad
L_{\si_0,r}\oplus L_{\si_1,r}\to W_{\si,r} \oplus K.
$$ 
(In the notation of [MS], $Q_{\si,r} = Q_{u_R,v_R}$.)  
As above there is a projection $pr_{\si_0,r}: L_{\si_0} \oplus L_{\si_1}
\to K$ and we write $Q_{\si,r}^Y$ for the restriction of $Q_{\si,r}$
to the fiber over $Y$.

As a guide to defining $Q_{a,app}$ consider the following diagram of spaces
 $$
\begin{array}{ccc} 
W_{\si,r}\oplus K&{\longleftarrow} & L_{\si_0,r}\oplus L_{\si_1,r}\\ 
\downarrow & &
\uparrow\\
 W_a\oplus K &\stackrel{Q_{a,app}}{\longleftarrow} &
L_a.\end{array} 
$$
where the maps are given by:
$$
\begin{array}{ccc} 
(\xi_0,\xi_1,Y)&\stackrel{Q_{\si,r}}{\longleftarrow} &(\eta_0,\eta_1) \\
\downarrow & & \uparrow\\ (\xi, Y)\in W_a\oplus K
&\stackrel{Q_{a,app}}{\longleftarrow} &\eta\in L_a.\end{array} 
$$
We define the horizontal arrow $Q_{a,app}$ by following the other three
arrows.   Here $\eta_\al$, for $\al = 0,1$, is the restriction of
$\eta$ to $\Si_{\si,\al} - D_{w}^\al(r)$ extended by $0$ as in [MS].  Note that
the $\eta_\al$ are in $L^p$ even  though they are not
continuous.  Next, decompose
$$
\begin{array}{ccl}
\eta_0 & = & \eta_0' + \io_{\si_0,r}(Y)\in \left(\im D_{\si_0,r}\right)
+ \io_{\si_0,r}(K) = L_{\si_0, r}\\
\eta_1 & = & \eta_1' + \io_{\si_1,r}(Y)  \in L_{\si_1,r}.
\end{array}
$$
Then $(\xi_0,\xi_1)  = Q_{\si,r}^Y(\eta_0',\eta_1')$. 
 Note that $\xi_0(w) =
\xi_1(w) = v$ say.   We then define the section $\xi$ by putting it equal
to $ \xi_\al$  on $\Si_\al - D_{w}^\al(r/\de)$  for $\al = 0,1$ and then extending it
over the neck using cutoff functions
so that it equals $\xi_0 + \xi_1 - v$ on the circle $\p D_w^0(r) = \p
D_w^1(r)\subset \Si_a$. 
In the formula below we think of the gluing map $\Psi_a$ of
Proposition~\ref{glue}  as inducing identifications
$$
\begin{array}{cccccc}
\Psi_a: \quad & A_0 & = &  D_{w}^0( r/\de) - D_{w}^0(r) &\to& D_{w}^1( r) -
D_{w}^1(r\de),\\
\Psi_a : \quad &
A_1 & = &  D_{w}^0( r) - D_{w}^0(r\de) &\to & D_{w}^1( r/\de) -
D_{w}^1(r).
\end{array}
$$
Then $\xi$ is given by 
$$
\xi(z) =\left\{ \begin{array}{ll} \xi_0(z) + (1 -
\be(z\de/r))(\xi_1(\Psi_a(z)) - v) & \mbox{if } z\in A_0,\\ 
\xi_1(\Psi_a(z)) + (1 - \be(z/r)(\xi_0(z) -
v) & \mbox{if } z\in A_1\end{array}\right.
$$
where $\be:\C\to [0,1]$ is a cutoff function that $=1$ if $|z|\le \de$ and
 $=0$ for $|z|\ge 1$.  

The next lemma is the analog of Lemma A.4.2 in [MS].  It shows that
$Q_{a,app}$ is the approximate inverse that we are seeking.
The norms used are the usual $L^p$-norms with respect to the
chosen metric on $\Si_\si$ and the glued metrics on $\Si_a$.
Note that we suppose that the metrics on $\Si_a$ agree with 
the standard model $\chi_r(|x|)|dx|^2$  on the annuli $D_{i}^\al(r/\de) - 
D_{i}^\al(r\de)$ (where $r= \|a\|^2$) so that $\Psi_a$ is an isometry.

\begin{lemma}\label{le:de} For all sufficiently small $\de$ there is
$r(\de)> 0$ and a cutoff function $\be$ such that for  all  $\eta \in L_a$,
$\|a\| \le r(\de)^2$, we have
 $$ \| (Dh_a\oplus \io_a) Q_{a,app} \eta -
\eta\| \le \frac 12 \|\eta\|. $$
\end{lemma}
\proof{}
It follows from the definitions that 
$$
(Dh_a\oplus \io_a) Q_{a,app} \eta =
\eta
$$
on each set $\Si_j - D_w^j(r/\de)$.  (Observe that $h_{\si_j,a} = h_a$ on this
domain so that $\io_{\si_j,a} = \io_a$ here.)
Therefore, we just have to consider
what happens on the subannuli 
$$
A_0 = D_w^0(r/\de) - D_w^0(r),\qquad \Psi_a(A_1) =
\Psi_a\left(D_w^0(r) - D_w^0(r\de)\right)
$$
of $\Si_a$.  In this region the maps $h_{\si_j,a}$ as well
as the glued map $h_a$ are constant so that the maps $\io_{\si_j,a}, \io_a$ are
constant. Further, the
linearizations $Dh_{\si_j,a}$ and $Dh_a$ are all equal  and on functions coincide
with the  usual  $\op$-operator.  We will consider what
happens in $\Psi_a(A_1)$, leaving the similar case of $A_0$ 
 to the reader.  

It is not hard to check that for $z\in A_1$
$$
\begin{array}{lclcc}
Dh_a(\xi_0(z)) & = & \eta_0'(z)  & = & -\io_a(Y),\\
Dh_a(\xi_1(\Psi_az))  & = &  \eta_1'(\Psi_a z). & &
\end{array}
$$
Let us write $\be_r$ for
the function  $
\be_r(z) = \be(z/r).
$
Then, if $r^2 = \|a\|$ and $(\xi,Y) = Q_{a,app}\eta$, we have for $z\in A_1$
\begin{eqnarray*}
(Dh_a \xi   + \io_a Y- \eta)(\Psi_a z) & = & \eta_1' (\Psi_a z)
+ (1 - \be_r)(-\io_a Y - Dh_a(v))(z) -\\& &
\qquad\quad  \op(\be_r)\otimes(\xi_0 - v)(z)
+ (\io_a Y -\eta)(\Psi_a z) \\
& = & (\be_r - 1) (\io_aY + Dh_a(v))(z) - \op(\be_r)\otimes(\xi_0 - v)(z).
\end{eqnarray*}
Therefore, taking the $L^p$-norm 
\begin{eqnarray*}
\|\left(Dh_a \xi   + \io_a Y- \eta\right)\circ \Psi_a\|_{L^p, A_1} & \le &
\|\io_a(Y)\|_{L^p,A_1} + \|Dh_a(v))\|_{L^p,A_1} +\\& &
\qquad\quad  \|\op(\be_r)\otimes(\xi_0 - v)\|_{L^p,A_1}.
\end{eqnarray*}
If $r$ is sufficiently small we can, by Lemma~\ref{le:NJ}
suppose that $\|\io_a(Y)\|_{L^p,A_1} \le \|\eta\|/12.$  Moreover,
because $v$ is a constant section $Dh_a$ acts on $v$ just by its zeroth
order part and so there are constants $c_1,c_2$ such that 
$$
\|Dh_a(v))\|_{L^p,A_1}
\le c_1\|v\| ({\rm area\,}A_1)^{1/p}\le c_2 \|v\| r^{2/p}.
$$
Furthermore, by [MS] Lemma A.1.2, given any $\eps > 0$ we can choose
$\de_\eps> 0$ and  $\be$ so that
$$
\|\op(\be_r)\otimes(\xi_0 - v)\|_{L^p,A_1} \le \eps \|\xi_0 -v\|_{W^{1,p}},
$$
for all $\de \le \de_\eps$.
Hence
\begin{eqnarray*}
\|Dh_a \xi   + \io_a Y- \eta\|_{L^p, A_1} & \le &
(c_2 r^{2/p} + \eps)(\|v\| + \|\xi_0 -v\|_{W^{1,p}}) + \|\eta\|/12\\
& \le & c_3(r^{2/p} + \eps)(\|(\eta_0', \eta_1')\|_{L^p}  + \|\eta\|/12\\
& \le & c_4(r^{2/p} + \eps)(\|(\eta_0, \eta_1)\|_{L^p}+ \|\eta\|/12\\
& = & \left(c_4(r^{2/p} + \eps) + 1/12\right)\|\eta\|,
\end{eqnarray*}
where the second inequality holds because of the uniform estimate for
the right inverse $Q_{\si,r}^Y$ and the 
third inequality holds because the projection of 
$L_{\si_0,r}\oplus L_{\si_1,r} $ onto the
subspace $ {\im} Dh_{\si,r} \oplus L_{\si_1,r}$ is continuous.
Then if we choose $\de_\eps$ so small that $c_4 \eps < 1/12$ and $r<<
\de_\eps$ so small that $ c_4 r < 1/12$ we find
$$
\|\left(Dh_a \xi   + \io_a Y- \eta\right)\circ\Psi_a\|_{L^p, A_1} \le 1/4 \|\eta\|.
$$
Repeating this for $A_0$ gives the desired result.\QED

Finally, we define the right inverse $Q_a$ by setting
$$
Q_a = Q_{a,app}\left((Dh_a\oplus \io_a) Q_{a,app}\right)^{-1}.
$$
It follows easily from the fact that the inverses $Q_{\si,r}$ are uniformly
bounded for $0 < r \le r_0$ that the $Q_a$ are too.

It remains to remark that the above construction can be carried out in such
as way as to be $\Ga_\si$-equivariant.  The only choice
left unspecified above was that of the right inverse  $Q_{\si,r}$.  This in
turn is determined by the choice of a subspace $R_{\si,r}$ of 
$$
 W^{1,p} (\Si, h_{\si,r}^*(TX))
$$
complementary to the kernel of $Dh_{\si,r}$.   But since $\Ga_\si$ is finite,
we can arrange that  $R_{\si,r}$ is $\Ga_\si$-equivariant.  For example,
since  $Dh_{\si,r}$ is a finite-dimensional space consisting of $C^\infty$
sections, we can take $R_{\si,r}$ to be the $L^2$-orthogonal complement
of  $Dh_{\si,r}$ defined with respect to a $\Ga_\si$-invariant norm on 
$h_{\si,r}^*(TX)$.\footnote{
As pointed out in [FO], the map
$$
\xi\mapsto \xi - \sum_j\langle \xi, e_j\rangle e_j, \quad \xi\in W_{\si,a}
$$
is well defined whenever $e_1,\dots, e_p$ is a finite set of
$C^\infty$-smooth sections.} 
Note that because $h_{\si,r} = h_{\si,r}\circ \ga$ for $\ga \in \Ga_\si,$ 
we can obtain a  $\Ga_\si$-invariant norm on 
$h_{\si,r}^*(TX)$  by integrating the pull-back by $h_{\si,r}$ of any norm on
the tangent bundle $TX$ with respect to a  $\Ga_\si$-invariant
area form on the domain $\Si_\si$.  We can achieve this uniformly over
$\Zz_J$  by choosing a suitable metric on each domain $\Si_\si$ as described at
the beginning of \S4.2.1. \MS

\begin{remark}\label{rmk:gen}\rm
If one is gluing  two branch components $\Si_{\ell_j}, j = 0,1$ of  $\Si$ then
both linearizations $Dh_{\ell_j}$ are surjective and one can construct the
inverse $Q_a$ to have image in $W_a$, thus forgetting about the summand
$K$.  The general gluing argument   combines both these cases.
If one is gluing at $N$ different points then one needs to choose $r$ so small
that one has an inequality of the form
$$
\|Dh_a \xi   + \io_a Y- \eta\|_{L^p, A} \le 1/4N \|\eta\|.
$$
on each of the $2N$-annuli $A$.  Note that the number of components of
$\Si$ is  bounded above  by some number that depends on $m$
(where $J\in \Jj_m$).  Hence there is always $r_0 > 0$  such that gluing at $\si$
is possible for all  $r < r_0$,  provided that one is looking a  family of
parametrized stable maps $\si = (\Si, h_\si, J)$ that is compact for each $J$
and where  $J\in\Jj_m$ is bounded in
$C^\infty$-norm.
\end{remark}

This completes the proof of Lemma~\ref{bound}  and hence of
Proposition~\ref{exun}. \MS

\NI
{\bf \S 4.2.4\, $\Aut^K(\Si)$-equivariance of $\TGg$}\MS

Note that there is an action of $S^1$ on the pair $(h_\si, a)$ that rotates 
one of the components (say $\Si_1$) of $\Si = \Si_0\cup \Si_1$ fixing the
intersection point  $w = \Si_0\cap \Si_1$.   We claim  that
by choosing an invariant metric on $\Si$ we can
make the whole construction invariant with respect to the action of this
compact group, i.e. so that  as {\it unparametrized} stable maps
$$
[\Si_a, \,\TGg(h_\si,a)] = [\Si_{\th\cdot a},\, \TGg( h_\si\cdot \th^{-1}, \th\cdot
a)].
$$
For then
there is an isometry $\psi$ from the glued domain $\Si_a$  to $\Si_{\th\cdot
a}$ such that 
$$
h_{a} = h^{\th}_{\th\cdot a}\circ\psi,
$$
where $h^{\th}_b $ denotes the pregluing of $h^\th_\si =
h_\si\cdot \th^{-1}$ with parameter $b$.   There is a similar formula for the
maps $h_{\si, r}$.  It is not hard to check that the rest of the construction can
be made compatible with this $S^1$ action.  It is important here to use the
Fukaya--Ono choice of $R_{\si, r}$ as described above, instead of cutting down
the domain of $Dh$  fixing the images of certain points as in ~[LiT], [LiuT1],
[Sieb].  

More generally, consider a parametrization $\tsi = (\Si, h)$ and an  arbitrary
element $\si
 \in
\Zz_J$.   Recall that a component of $\Si$ is said to be unstable if it contains
less than three special points, i.e. points where two
components of $\Si$ meet.  Each unstable branch component  has at least one
special point where it attaches to the rest of $\Si$ and so the identity
component of its automorphism group has the homotopy type of 
a circle.  Therefore, if there are $k$ such unstable components,
the torus group $T^k$ is a subgroup of $\Aut(\Si)$.
It is not hard to see that if the automorphism group $\Ga_\tsi$ of $\tsi$ is
nonzero, we can choose the action of $T^k$ to be $\Ga_\tsi$-equivariant
so that the groups fit together to form the compact group $\Aut^K(\Si)$ 
of Definition~\ref{def:gp}.

 Note further that if $(\Si', h')$ is obtained from
$\tsi = (\Si, h)$ by gluing, then  $\Aut^K(\Si')$ can be considered as a
subgroup of  $\Aut^K(\Si)$.  To see this, suppose for example that $\Si'$ is
obtained by gluing $\Si_i$ to $\Si_{j_i}$ with parameter $a$, and that both
these components have at most one other  special point.  Then
  we can choose metrics on $\Si_i\cup\Si_{j_i}$ that are invariant under an
$S^1$ action in each component and so that the glued metric on
$\Si_a$ is invariant under  the action of an $S^1$
in $\Aut^K(\Si')$.  Note that the diagonal subgroup $S^1\times S^1$
of $\Aut^K(\Si)$
acts trivially on the gluing parameters at the double point
$\Si_i\cap\Si_{j_i}$  since it rotates in opposite directions in the two tangent
spaces. It is now easy to check that if we write $\hat{\theta}$ for the image of
$\theta\in S^1$  in the diagonal subgroup of $S^1\times S^1$ then
 $$
(\Si_a, h'\circ\theta) = (\Si_a,\TGg(h,a)\circ\theta) = (\Si_a,
\TGg(h\circ\hat{\theta}, a)).
$$
Observe also that if $b$ is a gluing parameter at the intersection of $\Si_i$ with
some other component $\Si_k$ of $\tsi$, then it can also be considered as a
gluing parameter for $\tsi'$.  Moreover under this correspondence 
$\hat{\theta}\cdot b$ corresponds to $\theta\cdot b$.

These arguments prove the following result.

\begin{lemma}\label{le:repr}  Let $\tsi = (\Si, h)$  and suppose
that  the metric on $\Si$ is
$\Aut^K(\Si)$-invariant.  Then the following statements hold.\SS

\NI
(i) The composite $\Gg$ of $\TGg$ with the forgetful
map into the space of unparametrized stable maps is
$\Aut^K(\Si)$-invariant.  \SS

\NI
(ii) Divide the set $P$ of double points of $\Si$ into two sets $P_b, P_s$ 
and correspondingly write the gluing parameter $a$ as 
$ a_b + a_s$.  Suppose that
 $(\Si', h') = (\Si_{a_b}, \TGg(h,a_b))$, and consider $a_s$ as a gluing
parameter at $\tsi'$.   Then  one can choose metrics and choose the groups
$\Aut^K(\Si), \Aut^K(\Si')$ so that there is an inclusion
$$
\Aut^K(\Si')\to \Aut^K(\Si): \quad\theta\mapsto \hat{\theta}
$$ 
such that
$$
(\Si_{a_b}, h'\circ\theta^{-1}; \theta \cdot a_s) = 
(\Si_{a_b}, \TGg(h\circ\hat{\theta}^{-1}, a_b); \hat{\theta}\cdot a_s).
$$
Further, this can be done continuously as  $a_b$ (and hence $h'$) varies, and
smoothly if $\Si_{a_b}$ varies in a fixed stratum.  
\end{lemma}
\SS

\begin{remark}\rm  In their new paper~[LiuT2] \S5, Liu--Tian also  develop
a version of gluing that is invariant with respect to a partially defined torus
action. \end{remark}

 \NI
{\bf \S 4.2.5\, Globalization}\MS

The preceding paragraphs construct the gluing map $\TGg(h_\si,a)$  
 over a neighborhood $\Nn(\si)$ of one point $\si\in \Zz_J$.
 We now show how
to define a gluing map $\Gg_J:\Nn_\Vv(\Zz_J)\to \oMm(A-kF, \Kk_J)$ on 
a whole neighborhood  $\Nn_\Vv(\Zz_J)$ of $\Zz_J$ in the space of gluing
parameters $\Vv_J$.
  The only difficulty  in doing this lies in choosing a suitable
parametrized representative $s(\si) = (\Si, h_\si)$ of
the equivalence class $\si = [\Si, h]$ as $\si$ varies over $\Zz_J$.  
In other words,  in order to define  $\TGg(h_\si,a)$ we need to choose a
parametrization $h_\si: \Si \to X$ of the stable map $\si$, and
now we have to choose this consistently as $\si$ varies.  We now show 
that although we may not be able to make  a singlevalued choice $s(\si) =
h_\si$ continuously over $\Zz_J$ we can find a section that at each point  is
well defined modulo  the action of a suitable subgroup of  $\Aut^K(\Si)$.
More precisely, we claim the following.

\begin{lemma}\label{choice}  We may choose a continuous family of metrics 
$g_\si$ on $\Si_\si$ for $\si\in \Zz_J$ and
a family of parametrizations $s(\si)$ for each $\si\in \Zz_J$ such that

\NI
(i)
 $s(\si)$ consists of a $G_\si$-orbit of maps $h_\si:\Si_\si \to X$ 
and $g_\si$ is $G_\si$-invariant, where $G_\si\subset \Aut^K(\Si)$;

\NI
(ii)  the assignment $\si \to s(\si)$ is continuous in the sense that near
each $\si$ there is a  (singlevalued) continuous map $\si\to h_\si\in s(\si)$
whose restriction to each stratum is  smooth. Moreover, $g(\si)$
varies smoothly on each stratum.
\end{lemma} 
\proof{}
The strata in $\Zz_J$  can be partially ordered with $\Ss' \le \Ss$ if there
is a  gluing that takes an element in the stratum $\Ss$ to that in
$\Ss'$, i.e. if the stratum $\Ss$ is contained in the closure of $\Ss'$.   If $\Ss$
is maximal under this ordering and $\si \in \Ss$, then each branch
component in $\Si$ is mapped to a fiber by a map of degree $\le 1$.  It is
easy to check that in this case there is a unique identification  of the  domain
$\Si_\si$ with a union of spheres such  that the map $h_\si$ is either
constant or is the identity map on each branch component and a section on
the stem: cf. Example~\ref{ex:lens}.   We assume this done, and then extend
the choice of parametrization to a neighborhood of each of these maximal
strata by gluing.  

We now start extending our choice $s(\si) = h_\si$ of
parametrization  to the whole of $\Zz_J$ 
 by downwards induction over the partially ordered strata. Clearly we can
always choose a parametrization modulo the action of $\Aut^K(\Si)$.  
In order for the image of the fiber $\pi_J^{-1}(\si) = \{(\si,a): |a|<\eps\}$
 under  the gluing map to be
independent of this choice, we need the metric $g_\si$ on $\Si_\si$ to be
$\Aut^K(\Si_\si)$-invariant.  This choice of metric can be assumed to be smooth
as $\si$ varies in a stratum.  However, it 
cannot always be chosen continuously as $\si$ goes
from one stratum to another.  For example, if $\si$ has one component
$\Si_i$ with $3$ special points at $0,1,\infty$ and that is glued to
some component $\Si_{j_i}$ at $1$ with gluing parameter $a$, then the
resulting component $\Si_a$ is unstable if $\Si_{j_i}$ has no other special
points.  But for small $|a|$ the metric on $\Si_a$  is determined by the metrics
on $\Si = \Si_i\cup\Si_{j_i}$ by the gluing construction and cannot be chosen
to be $S^1$-invariant.
On the other hand, if  both $\Si_i$ and $\Si_{j_i}$ have at most one
other  special point then the glued metric on $\Si_a$ will be
$S^1$-invariant provided that the original metrics on $\Si_i, \Si_{j_i}$ were
also $S^1$-invariant.

The above remarks show that suitable $g_\si^\Ss, s(\si)^\Ss$ and $G_\si^\Ss$
can be defined over each stratum $\Ss$, and in particular over 
maximal strata.  If  $g_\si, s(\si) $ and $G_\si$
 are already suitably defined over  
some union $Y$ of strata, then the above 
remarks about gluing  show that they can
be extended to a neighborhood $\Uu(Y)$ of $Y$.
Let us write $g_\si^{gl}, s(\si)^{gl}$ and $G_\si^{gl}$ for
the objects obtained by gluing when $\si\in \Uu(Y)$.
Then, if $\be: \Uu(Y)\cup \Ss \to [0,1]$ is a smooth cutoff function
that equals $0$ near $Y$ and $1$ near the boundary of $\Uu(Y)$, set
\begin{eqnarray*}
g_\si & = &  (1-\be(\si))g_\si^{gl} + \be(\si) g_\si^\Ss, \si \in \Ss\\
s(\si) & =& s(\si)^\Ss,\quad \mbox{if }\;\be(\si)=1,\\
& = & s(\si)^{gl}\quad\mbox{otherwise},\\
G_\si & =& G_\si^\Ss,\quad \mbox{if }\;\be(\si)=1,\\
& = & G_\si^{gl}\quad\mbox{otherwise}.
\end{eqnarray*}
It is easy to check that the required conditions are satisfied.\QED

\MS

\NI{\bf Proof of Proposition~\ref{glue2}}
\MS

By Lemmas~\ref{le:repr} and~\ref{choice} there is a well
defined continuous gluing map
$$
\Gg_J: \Nn_\Vv(\Zz_J) \to \oMm(A-kF, \Kk_J).
$$
that  restricts
on $\Zz_J$ to the inclusion.  Therefore, because $\Zz_J$ is compact, the
injectivity of $\Gg_J$ on a small neighborhood $\Nn_\Vv(\Zz_J)$ follows from
the  local injectivity statement in 
Proposition~\ref{exun}.  Similarly, the local
surjectivity of Proposition~\ref{exun} 
implies that the image of $\Gg_J$ is open
in $\oMm(A-kF,\Kk_J)$.
Note that all the restrictions made on the size of $\Nn_\Vv(\Zz_J)$ vary
smoothly with $J$ (and involve no more than the $C^2$ norm of $J$).  Hence 
$\cup_J\im \Gg_J$ is an open subset of 
$\oMm(A-kF,\Jj)$. \QED 
\MS\MS

\NI{\bf References}\MS

\NI
[A]  M. Abreu,  Topology of symplectomorphism groups of $S^2\times S^2$,
{\it Invent. Math.} {\bf 131}, (1998), 1--23.

\MS

\NI
[AM]  M. Abreu and D. McDuff,    Topology of symplectomorphism groups
of rational ruled surfaces, in preparation.\MS

\NI
[FO]   K. Fukaya and K. Ono, Arnold conjecture and
Gromov--Witten invariants, to appear in {\it Topology}

\MS

\NI [HLS] H.~Hofer, V.~Lizan and J.-C.~Sikorav,   On genericity for complex
curves in $4$-dimensional, almost complex manifolds,
\MS

\NI
[HS] H. Hofer and D. Salamon,
Gromov compactness and stable maps, 
Preprint (1997).
\MS

\NI
[K] P. Kronheimer, Some nontrivial families of symplectic structures,
preprint (1998)
\MS

\NI [LM]   F. Lalonde and D McDuff, $J$\/-curves and the 
       classification of rational and ruled
        symplectic $4$\/-manifolds,  in {\it Contact and Symplectic
Geometry}, ed C. Thomas,
       Camb Univ Press (1996).

\MS

\NI
[LiT]  Jun Li and Gang Tian, Virtual moduli cycles and Gromov--Witten
invariants of general symplectic manifolds, preprint (1996)
\MS

\NI
[LiuT1] Gang Liu and Gang Tian, Floer homology and Arnold
conjecture, to appear in {\it Journ. Differential Geometry}\MS

\NI [LiuT2] Gang Liu and Gang Tian, Weinstein conjecture and GW invariants,
preprint (1997).\MS

 \NI
[Lo] W. Lorek, Generalized Cauchy--Riemann operators in Symplectic
Geometry, Ph. D. thesis, Stony Brook (1996).\MS

 \NI [MP]
D. McDuff and L. Polterovich,  
Symplectic packings and algebraic geometry,
{\it Inventiones Mathematicae}, {\bf 115}, (1994) 405--29.

\NI 
[MS]   D. McDuff and D.A. Salamon, {\it $J$-holomorphic curves and
quantum cohomology}, Amer Math Soc Lecture Notes \#6, Amer.
Math. Soc. Providence (1995).
\MS

\NI
[R]  Yongbin Ruan,  Virtual neighborhoods and pseudoholomorphic 
curves, preprint
Preprint alg-geom/9611021

\MS

\NI
[S]  B. Siebert,  Gromov--Witten invariants for general
symplectic manifolds. Preprint, Bochum (1996).
\end{document}